\documentstyle[12pt]{article}

  \input amssym.def \input amssym \def\mb{\mbox}
  \newfont{\bbb}{msbm10 scaled\magstep1}
  \newfont{\blb}{msbm7 scaled\magstep1}

\def\baselinestretch{1} \topmargin -12pt \headsep 0pt \footskip 30pt
\def\medskipamount{12pt}

\newcommand{\bit}[1]{\section{#1} \setcounter{equation}{0}}
\renewcommand{\theequation}{\thesection .\arabic{equation}}

\newcommand{\re}[1]{\mbox{\rm \bf (\ref{#1})}}

\catcode`\@=\active
\catcode`\@=11
\def\@eqnnum{\hbox to .01pt{}\rlap{\bf \hskip -\displaywidth(\theequation)}}
\catcode`\@=12

\newenvironment{s}[1]
{ \addvspace{\medskipamount} \refstepcounter{equation}
\noindent {\bf (\theequation) #1.} \begin{em}}
{\end{em} \par \addvspace{\medskipamount} }

\newenvironment{r}[1]
{ \addvspace{\medskipamount} \refstepcounter{equation}
\noindent {\bf (\theequation) #1.} }
{\par \addvspace{\medskipamount} }

\begin{document}


\catcode`\@=\active
\catcode`\@=11
\newcommand{\nc}{\newcommand}


\nc{\vars}[2]
{{\mathchoice{\mb{#1}}{\mb{#1}}{\mb{#2}}{\mb{#2}}}}
\nc{\Aff}{\vars{\bbb A}{\blb A}}
\nc{\C}{\vars{\bbb C}{\blb C}}
\nc{\Hyp}{\vars{\bbb H}{\blb H}}
\nc{\N}{\vars{\bbb N}{\blb N}}
\nc{\Pj}{\vars{\bbb P}{\blb P}}
\nc{\Q}{\vars{\bbb Q}{\blb Q}}
\nc{\R}{\vars{\bbb R}{\blb R}}
\nc{\Z}{\vars{\bbb Z}{\blb Z}}


\nc{\oper}[1]{\mathop{\mathchoice{\mbox{\rm #1}}{\mbox{\rm #1}}
{\mbox{\scriptsize \rm #1}}{\mbox{\tiny \rm #1}}}\nolimits}
\nc{\Aut}{\oper{Aut}}
\nc{\chr}{\oper{char}}
\nc{\diag}{\oper{diag}}
\nc{\End}{\oper{End}}
\nc{\Fl}{\oper{Fl}}
\nc{\Gr}{\oper{Gr}}
\nc{\Hom}{\oper{Hom}}
\nc{\NS}{\oper{NS}}
\nc{\Par}{\oper{Par}}
\nc{\Pic}{\oper{Pic}}
\nc{\Proj}{\oper{Proj}}
\nc{\Quot}{\oper{Quot}}
\nc{\Spec}{\oper{Spec}}
\nc{\tr}{\oper{tr}}


\nc{\GL}[1]{{\rm GL(#1)}}
\nc{\PSL}[1]{{\rm PSL(#1)}}
\nc{\SL}[1]{{\rm SL(#1)}}


\nc{\ep}{\epsilon}
\nc{\ga}{\gamma}
\nc{\la}{\lambda}
\nc{\La}{\Lambda}
\nc{\si}{\sigma}


\nc{\A}{A}
\nc{\aff}{{k}}
\nc{\blowup}{\tilde{\Aff}^{\raisebox{-.5ex}{$\scriptstyle 2$}}}
\nc{\down}{\Big\downarrow}
\nc{\be}{{\bf E}}
\nc{\beqas}{\begin{eqnarray*}}
\nc{\beqa}{\begin{eqnarray}}
\nc{\beq}{\begin{equation}}
\nc{\bl}{\vskip 1.2ex }
\nc{\ci}{{\Im}}
\nc{\co}{{\cal O}}
\nc{\E}{E}
\nc{\eeqas}{\end{eqnarray*}}
\nc{\eeqa}{\end{eqnarray}}
\nc{\eeq}{\end{equation}}
\nc{\emb}{\hookrightarrow}
\nc{\fp}{\mbox{     $\Box$}}
\nc{\half}{\frac{\scriptstyle 1}{\scriptstyle 2}}
\nc{\k}{{k}}
\nc{\kst}{{k^{\times}}}
\nc{\lrow}{\longrightarrow}
\nc{\m}{{\frak m}}
\nc{\mod}{/ \! \! /}
\nc{\pf}{{\em Proof}}
\nc{\sans}{\backslash}
\nc{\st}{\, | \,}

\catcode`\@=12



\noindent
{\LARGE \bf Geometric invariant theory and flips}
\medskip \\
{\bf Michael Thaddeus }\\
Mathematical Institute, 24--29 St Giles, Oxford OX1 3LB, England
\medskip
\smallskip


\noindent Ever since the invention of geometric invariant theory, it
has been understood that the quotient it constructs is not entirely
canonical, but depends on a choice: the choice of a {\em
linearization} of the group action.  However, the founders of the
subject never made a systematic study of this dependence.  In light of
its fundamental and elementary nature, this is a rather surprising
gap, and this paper will attempt to fill it.

In a way, this neglect is understandable, because the different
quotients must be related by birational transformations, whose
structure in higher dimensions is poorly understood.  However, it has
been considerably clarified in the last dozen years with the advent of
Mori theory.  In particular, the birational transformations that Mori
called {\em flips} are ubiquitous in geometric invariant theory;
indeed, one of our main results \re{1f} describes the mild conditions
under which the transformation between two quotients is given by a
flip.  This paper will not use any of the deep results of Mori theory,
but the notion of a flip is central to it.

The definition of a flip does not describe the birational
transformation explicitly, but in the general case there is not much
more to say.  So to obtain more concrete results, hypotheses must be
imposed which, though fairly strong, still include many interesting
examples.  The heart of the paper, \S\S4 and 5, is devoted to
describing the birational transformations between quotients as
explicitly as possible under these hypotheses.  It turns out that
there are fairly explicit smooth loci in two different quotients whose
blow-ups are isomorphic.  Thus the two quotients are related by a
blow-up followed by a blow-down.  This is somewhat at odds with the
point of view of Mori theory, which views a flip as two contractions,
not two blow-ups; but it facilitates explicit calculations of such
things as topological cohomology or Hilbert polynomials.

The last three sections of the paper put this theory into practice,
using it to study moduli spaces of points on the line, parabolic
bundles on curves, and Bradlow pairs.  An important theme is that the
structure of each individual quotient is illuminated by understanding
the structure of the whole family.  So even if there is one especially
natural linearization, the problem is still interesting.  Indeed, even
if the linearization is unique, useful results can be produced by
enlarging the variety on which the group acts, so as to create more
linearizations.

I believe that this problem is essentially elementary in nature, and I
have striven to solve it using a minimum of technical machinery.  For
example, stability and semistability are distinguished as little as
possible.  Moreover, transcendental methods, choosing a maximal torus,
and invoking the numerical criterion are completely avoided.  The only
technical tool relied on heavily is the marvelous Luna slice theorem
\cite{luna}.  Luckily, this is not too difficult itself, and there is
a good exposition in GIT, appendix 1D.  This theorem is used, for
example, to give a new, easy proof of the Bialynicki-Birula
decomposition theorem \re{2e}.

Section 1 treats the simplest case: that of an affine variety $X$
acted on by the multiplicative group $\kst$, and linearized on the
trivial bundle.  This case has already been treated by Brion and
Procesi \cite{bp}, but the approach here is somewhat different,
utilizing $\Z$-graded rings.  The four main results are models for
what comes later.  The first result, \re{2i}, asserts that the two
quotients $X \mod \pm$ coming from nontrivial linearizations are
typically related by a flip over the affine quotient $X \mod 0$.  The
second, \re{2d}, describes how to blow up ideal sheaves on $X \mod
\pm$ to obtain varieties which are both isomorphic to the same
component of the fibred product ${X \mod -} \times_{X \mod 0} {X \mod
+}$.  In other words, it shows how to get from $X \mod +$ to $X \mod
-$ by performing a blow-up followed by a blow-down.  The third result,
\re{2j}, asserts that when $X$ is smooth, the blow-up loci are
supported on subvarieties isomorphic to weighted projective fibrations
over the fixed-point set.  Finally, the fourth, \re{2k}, identifies
these fibrations, in what gauge theorists would call the quasi-free
case, as the projectivizations of weight subbundles of the normal
bundle to the fixed-point set.  Moreover, the blow-ups are just the
familiar blow-ups of smooth varieties along smooth subvarieties.

Sections 3, 4, and 5 are concerned with generalizing these results in
three ways.  First, the variety $X$ may be any quasi-projective
variety, projective over an affine.  Second, the group acting may be
any reductive algebraic group.  Third, the linearization may be
arbitrary.  But \S2 assumes $X$ is normal and projective, and is
something of a digression.  It starts off by introducing a notion of
$G$-algebraic equivalence, and shows, following Mumford, that
linearizations equivalent in this way give the same quotients.  Hence
quotients are really parametrized by the space of equivalence classes,
the $G$-N\'eron-Severi group $\NS^G$.  Just as in the
Duistermaat-Heckman theory in symplectic geometry, it turns out
\re{1b} that $\NS^G \otimes \Q$ is divided into chambers, on which the
quotient is constant.

The analogues of the four main results of \S1 then apply to quotients
in adjacent chambers, though they are stated in a somewhat more
general setting.  The first two results are readily generalized to
\re{1f} and \re{1g}.  The second two, however, require the hypotheses
mentioned above; indeed, there are two analogues of each.  The first,
\re{2n} and \re{1x}, make fairly strong hypotheses, and show that the
weighted projective fibrations are locally trivial.  The second,
\re{1l} and \re{1n}, relax the hypotheses somewhat, but conclude only
that the fibrations are locally trivial in the \'etale topology.
Counterexamples \re{1r} and \re{2p} show that the hypotheses are
necessary.

The strategy for proving all four of these results is not to imitate
the proofs in the simple case, but rather to reduce to this case by
means of a trick.  In fact, given a variety $X$ acted on by a group
$G$, and a family of linearizations parametrized by $t$, we construct
\re{1h} a new variety $Z$, dubbed the ``master space'' by Bertram, acted
on by a torus $T$, and a family of linearizations on $\co(1)$
parametrized by $t$, such that $X \mod G(t) = Z \mod T(t)$ naturally.
This reduces everything to the simple case.

The final sections, \S\S6, 7, and 8, are in a more discursive style;
they explain how to apply the theory of the preceding sections to some
examples.  In all cases, the strongest hypotheses are satisfied, so
the best result \re{1x} holds.  Perhaps the simplest interesting
moduli problems are those of (ordered or unordered) sets of $n$ points
in $\Pj^1$; these are studied in \S6.  The ideas here should have many
applications, but only a very simple one is given: the formula of
Kirwan \cite{k} for the Betti numbers of the moduli spaces when $n$ is
odd.  In \S7 the theory is applied to parabolic bundles on a curve,
and the results of Boden and Hu \cite{bhu} are recovered and extended.
Finally, in \S8, the theory is applied to Bradlow pairs on a curve,
recovering the results of the author \cite{t1} and Bertram et al.\
\cite{bdw}.

While carrying out this research, I was aware of the parallel work of
Dolgachev and Hu, and I received their preprint \cite{dh} while this
paper was being written.  Their main result is contained in the third
of the four main results I describe, \re{1l}; and of course, some of
the preliminary material, corresponding to my \S2, overlaps.  I am
indebted to them for the observation quoted just after \re{1y}, and
for the result \re{1s}, though my proof of the latter is original.
Dolgachev and Hu do not, however, include the results on flips or
blow-ups, study the local triviality of the exceptional loci, or
identify the projective bundles in the quasi-free case.  For them,
this is not necessary, since they appear \cite{hu} to be interested
chiefly in computing Betti numbers and intersection Betti numbers of
quotients, and for this, their main result suffices, together with the
deep results of Beilinson, Bernstein, and Deligne \cite{bbd}.  I am
more interested in computing algebraic cohomology, as in \cite{t1};
for this, a precise characterization of the birational map between
quotients is necessary, which led me to the present paper.  In any
case, even where our results coincide, our methods of proof are quite
different.

A few words on notation and terminology.  Many of the statements
involve the symbol $\pm$.  This should always be read as two distinct
statements: that is, $X^\pm$ means $X^+$ (resp.\ $X^-$), never $X^+
\cup X^-$ or $X^+ \cap X^-$.  Similarly, $X^\mp$ means $X^-$ (resp.\
$X^+$). The quotient of $X$ by $G$ is denoted $X \mod G$, or $X \mod
G(L)$ to emphasize the choice of a linearization $L$.  When there is
no possibility of confusion, we indulge in such abuses of notation as
$X \mod \pm$, which are explained in the text.  For stable and
semistable sets, the more modern definitions of \cite{n} are followed,
not those of \cite{mf}, which incidentally is often referred to as
GIT.  Points are assumed to be closed unless otherwise stated.  The
stabilizer in $G$ of a point $x \in X$ is denoted $G_x$.  When $E$ and
$F$ are varieties with morphisms to $G$, then $E \times_G F$ denotes
the fibred product; but if $G$ is a group acting on $E$ and $F$, then
$E \times_G F$ denotes the twisted quotient $(E \times F) / G$.
Unfortunately, both notations are completely standard.

All varieties are over a fixed algebraically closed field $\k$.  This
may have any characteristic: although we use the Luna slice theorem,
which is usually said \cite{luna,mf} to apply only to characteristic
zero, in fact this hypothesis is used only to show that the stabilizer
must be linearly reductive.  Since all the stabilizers we encounter
will be reduced subgroups of the multiplicative group, this will be
true in any characteristic.

By the way, most of the results in \S\S1 and 3 apply not only to
varieties, but to schemes of finite type over $\k$.  But for
simplicity everything is stated for varieties.

Finally, since the experts do not entirely agree on the definition of
a flip, here is what we shall use.  Let $X_- \to X_0$ be a {\em small
contraction} of varieties over $\k$.  This means a small birational
proper morphism; {\em small} means that the exceptional set has
codimension greater than $1$.  (This appears to be the prevailing
terminology in Mori theory \cite[(2.1.6)]{ko}, but it is called {\em
semismall\/} in intersection homology, where {\em small} has a
stronger meaning.)  Let $D$ be a $\Q$-Cartier divisor class on $X_-$
which is relatively negative over $X_0$; that is, $\co(-D)$ is
relatively ample.  Then the {\em $D$-flip} is a variety $X_+$, with a
small contraction $X_+ \to X_0$, such that, if $g: X_- \dasharrow X_+$
is the induced birational map, then the divisor class $g_*D$ is
$\Q$-Cartier, and $\co(D)$ is relatively ample over $X_0$.  We
emphasize the shift between ampleness of $\co(-D)$ and that of
$\co(D)$.  If a flip exists it is easily seen to be unique.  Note that
several authors, including Mori \cite{m}, require that each
contraction reduce the Picard number by exactly $1$.  We will not
require this; indeed, it is not generally true of our flips
\cite[4.7]{toric}.  For convenience, $D$-flips will be referred to
simply as {\em flips}.  However, in the literature, the unmodified
word {\em flip} has traditionally denoted a $K$-flip where $K$ is the
canonical divisor of $X_-$; this is not what we will mean.

\bit{The simplest case}

We begin by examining the simplest case, that of an affine variety
acted on by the multiplicative group of $\k$.  This has been studied
before in several papers, that of Brion and Procesi \cite{bp} being
closest to our treatment; but we will clarify, extend, and slightly
correct the existing results.

Let $R$ be a finitely-generated integral algebra over the
algebraically closed field $\k$, so that $X = \Spec R$ is an affine
variety over $\k$.  In this section only, $G$ will denote the
multiplicative group of $\k$.  An action of $G$ on $\Spec R$ is
equivalent to a $\Z$-grading of $R$ over $\k$, say $R = \bigoplus_{i
\in \Z} R_i$.  We will study geometric invariant theory quotients $X
\mod G$, linearized on the trivial bundle.

So choose any $n \in \Z$, and define a $\Z$-grading on $R[z]$ by $R_i
\subset R[z]_i$, $z \in R[z]_{-n}$.  Of course $R[z]$ is also
$\N$-graded by the degree in $z$, but this $\Z$-grading is different.
Since $X = \Spec R = \Proj R[z]$, the $\Z$-grading is equivalent to a
linearization on $\co$ of the $G$-action on $X$.  The quotient is
$\Proj R[z]^{G} = \Proj R[z]_0 = \Proj \bigoplus_{i \in \N} R_{ni}
z^i$.  For $n = 0$, this is just $\Proj R_0[z] = \Spec R_0$, the usual
affine quotient \cite[3.5; GIT Thm.\ 1.1]{n}.  For $n > 0$, $\Proj
\bigoplus_{i \in \N} R_{ni} = \Proj \bigoplus_{i \in \N} R_{i}$ by
\cite[II Ex.\ 5.13]{h} (the hypothesis there is not needed for the
first statement); similarly for $n < 0$, $\Proj \bigoplus_{i \in \N}
R_{ni} = \Proj \bigoplus_{i \in \N} R_{-i}$.  Hence we need concern
ourselves only with the quotients when $n = 0$, $1$, and $-1$; we
shall refer to them in this section as $X \mod 0$, $X \mod +$ and $X
\mod -$ respectively. Note that $X \mod \pm$ are projective over $X
\mod 0$.

\begin{s}{Proposition}
\label{2h}
If $X \mod + \neq \varnothing \neq X \mod -$, then the natural
morphisms $X \mod \pm \to X \mod 0$ are birational.
\end{s}

Of course, if say $X \mod - = \varnothing$, then $X \mod +$ can be any
$\k$-variety projective over $X \mod 0$.

\pf.  For some $d > 0$ $R_{-d}$ contains a nonzero element $t$.  The
function field of $X \mod 0$ is $\{ r/s \st r,s \in R_0 \}$, while
that of $X \mod +$ is $\{ r/s \st r,s \in R_{di} \mbox{ \rm for some }
i \geq 0 \}$.  But the map $r/s \mapsto (rt^i)/(st^i)$ from the latter
to the former is an isomorphism.  \fp \bl

Let $X^\pm \subset X$ be the subvarieties corresponding to the ideals
$\langle R_i \st \mp i > 0 \rangle$ (note the change of sign), and let
$X^0 = X^+ \cap X^-$; then $X^0$ corresponds to the ideal $\langle R_i
\st i \neq 0 \rangle$.  Also say $\lim g \cdot x = y$ if the morphism
$G \to X$ given by $g \mapsto g \cdot x$ extends to a morphism $\Aff^1
\to X$ such that $0 \mapsto y$.

\begin{s}{Lemma}
\label{2l}
As sets, $X^\pm = \{ x \in X \st \exists \lim g^{\pm 1} \cdot x \}$,
and $X^0$ is the fixed-point set for the $G$-action.
\end{s}

\pf.  A point $x$ is in $X^+$ if and only if for all $n < 0$, $R_n$ is
killed by the homomorphism $R \to \k[x, x^{-1}]$ of graded rings
induced by $g \mapsto g \cdot x$.  This in turn holds if the image of
$R$ is contained in $\k[x]$, that is, if $G \to X$ extends to $\Aff^1
\to X$.  The proof for $X^-$ is similar.

Hence $x \in X^0$ if and only if $\lim g^{\pm 1} \cdot x$ both exist,
that is, the closure of $G \cdot x$ is a projective variety in $X$.
Since $X$ is affine, this means $x$ is a fixed point of $G$. \fp

\begin{s}{Proposition}
\label{2a}
\mbox{\rm (a)} $X^{ss}(0) = X$; \mbox{\rm (b)} $X^s(0) = X \sans (X^+
\cup X^-)$ ; \mbox{\rm (c)} $X^{ss}(\pm) = X^s(\pm) = X \sans X^\mp$.
\end{s}

\pf.  Recall that $x \in X^{ss}(L)$ if for some $n > 0$ there exists
$s \in H^0(L^n)^G$ such that $s(x) \neq 0$, and $x \in X^{s}(L)$ if
the morphism $G \to X^{ss}(L)$ given by $g \mapsto g \cdot x$ is
proper.  For $L = 0$, $H^0(L^n)^{G} = R_0$ for all $n$, but this
contains 1, which is nowhere vanishing.  That is all there is to (a).
The valuative criterion implies that the morphism $G \to X$ is proper
if and only if the limits do not exist, which together with \re{2l}
implies (b).  For $L = \pm$, $H^0(L^n)^{G} = R_{\pm n}$, so
$X^{ss}(\pm) = X \sans X^\mp$ follows immediately from the definition
of $X^\pm$.  The additional condition of properness needed for $x \in
X^s(\pm)$ is equivalent, by the valuative criterion, to $\lim g \cdot
x$ and $\lim g^{-1} \cdot x \notin X^{ss}(\pm)$. But one does not
exist, and the other, if it exists, is fixed by $G$, so is certainly
not in $X^{ss}(\pm)$.  \fp

\begin{s}{Corollary}
\label{2m}
The morphisms $X \mod \pm \to X \mod 0$ are isomorphisms on the
complements of $X^{\pm} \mod \pm \to X^\pm \mod 0$.
\end{s}

In good cases, $X^{\pm} \mod \pm$ will be exactly the exceptional loci
of the morphisms, but they can be smaller, even empty---for instance
$X \mod - \to X \mod 0$ in \re{2g} below.

\pf.  By \re{2a}, the sets $(X \sans X^\pm) \mod \pm$ and $(X \sans
X^\pm) \mod 0$ contain no quotients of strictly semistable points.
They are therefore isomorphic.  \fp

\begin{s}{Proposition}
There are canonical dominant morphisms $\pi_\pm: X^\pm \to X^0$ such
that for all $x \in X^\pm$, $\pi_\pm (x) = \lim g^{\pm 1} \cdot x$.
\end{s}

\pf.  Note first that $R_0 \cap \langle R_i \st \pm i > 0 \rangle =
R_0 \cap \langle R_i \st i \neq 0 \rangle$.  So $R / \langle R_i \st i
\neq 0 \rangle = R_0 / \langle R_i \st i \neq 0 \rangle$ are naturally
included in $R / \langle R_i \st \pm i > 0 \rangle$ as the
$G$-invariant parts.  Hence there are natural dominant morphisms
$\pi_\pm: X^\pm \to X^0$.

Because $\pi_\pm$ is induced by the inclusion of the degree 0 part in
$R / \langle R_i \st \mp i > 0 \rangle$, $\pi_\pm (x)$ is the unique
fixed point such that $f(\pi_\pm (x)) = f(x)$ for all $f \in R_0$.
But $f \in R_0$ means it is $G$-invariant, hence constant on orbits,
so the same property is satisfied by $\lim g^{\pm 1} \cdot x$, which
is a fixed point in the closure of $G \cdot x$.  \fp \bl

The next two results digress to show what the results so far have to
do with flips.

\begin{s}{Proposition}
\label{2i}
If $X^{\pm} \subset X$ have codimension $\geq 2$, then the birational
map $f: X \mod - \dasharrow X \mod +$ is a flip with respect to
$\co(1)$.
\end{s}

\pf.  The hypothesis implies that the open sets $(X \sans X^\pm) \mod
\pm$ in $X \mod \pm$ have complements of codimension $\geq 2$.  But by
\re{2m} these open sets are identified by $f$.  Hence there is a
well-defined push-forward $f_*$ of divisors.

For some $n > 0$ the twisting sheaves $\co(\pm n) \to X \mod \pm$ are
line bundles.  Indeed, they are the descents \cite{dn} from $X$ to $X
\mod \pm$ of the trivial bundle $\co$, with linearization given by $n$
as at the beginning of this section.  Consequently, they agree on the
open sets $(X \sans X^\pm) \mod \pm$, so $f_* \co(-n) = \co(n)$.  But
$\co(\pm n) \to X \mod \pm$ are relatively ample over $X \mod 0$, so $f$
is a flip.  \fp

\begin{s}{Proposition}
\label{2o}
Let $Y_0$ be normal and affine over $\k$, and let $f: Y_- \dasharrow
Y_+$ be a flip of normal varieties over $Y_0$.  Then there exists $X$
affine over $\k$ and a $G$-action on $X$ so that $Y_0 = X \mod 0$,
$Y_\pm = X \mod \pm$.
\end{s}

\pf.  Let $L = \co(D)$, where $D$ is as in the definition of a flip.
Since $Y_\pm$ are normal, and the exceptional loci of $f$ have
codimension $\geq 2$, $f$ induces an isomorphism $H^0(Y_-, L^n) \cong
H^0(Y_+, f_*L^n)$ for all $n$.  The $\N$-graded algebras $\bigoplus_{n
< 0} H^0(Y_-, L^n)$ and $\bigoplus_{n > 0} H^0(Y_+, f_*L^n)$ are the
homogeneous coordinate rings of the quasi-projective varieties $Y_-$
and $Y_+$ respectively, so are finitely-generated over $\k$.  Hence
the same is true of the $\Z$-graded algebra $R = \bigoplus_{n \in \Z}
H^0(Y_-, L^n)$.  Let $X = \Spec R$ with the $G$-action coming from the
grading.  Then $X \mod 0 = Y_0$ and $X \mod \pm = Y_\pm$.  I thank
Miles Reid for pointing out this simple proof.  \fp \bl

In order to describe the birational map $X \mod - \dasharrow X \mod +$
more explicitly, we will next construct a variety birational to $X
\mod \pm$ which dominates them both.  It is admittedly true in general
that any birational map can be factored into a blow-up and blow-down
of some sheaves of ideals.  The virtue of the present situation,
however, is that these sheaves can be identified fairly explicitly,
and that the common blow-up is precisely the fibred product.

Choose $d > 0$ such that $\bigoplus_{i \in \Z} R_{di}$ is generated by
$R_0$ and $R_{\pm d}$.  Then let $\ci^\pm$ be the sheaves of ideals on
$X$ corresponding to $\langle R_{\mp d} \rangle$.  Let $\ci^\pm \mod
\pm$ on $X \mod \pm$ and $\ci^\pm \mod 0$ on $X \mod 0$ be the ideal
sheaves of invariants of $\ci^\pm$, that is, the sheaves of ideals
locally generated by the invariant elements of $\ci^\pm$.  Note that
$\ci^\pm$ are supported on $X^\pm$, so that $\ci^\pm \mod \pm$ are
supported on $X^\pm \mod \pm$.

For $i, j \geq 0$, let $R_{i,j} = R_i \cdot R_{-j} \subset
R_{i-j}$.

\begin{s}{Lemma}
\label{2b}
The ideal sheaf $(\ci^+ + \ci^-) \mod 0$ is exactly $\langle R_{d,d}
\rangle$, and its pullbacks by the morphisms $X \mod \pm \to X \mod 0$
are $\ci^\pm \mod \pm$.
\end{s}

\pf. The ideals in $R$ corresponding to $\ci^\pm$ are by definition
$\langle R_{\mp d} \rangle$, and $\langle R_{\mp d} \rangle \cap R_0 =
R_{d,d}$, so $(\ci^+ + \ci^-) \mod 0 = \langle R_{d,d} \rangle$.
Regard $X \mod \pm$ as quotients with respect to the linearizations
$\pm d$.  Then $X \mod \pm = \bigoplus_{i \geq 0} R_{\pm di}$, so for
any $\si \in R_{\pm d}$, $\Spec (\si^{-1} \bigoplus_{i \geq 0} R_{\pm
di})_0$ is an affine in $X \mod \pm$.  But $\ci^\pm \cap (\si^{-1}
\bigoplus_{i \geq 0} R_{\pm di})_0 = \si^{-1} \langle R_{d,d}
\rangle$, so locally $\ci^\pm \mod \pm$ is the pullback of $\langle
R_{d,d} \rangle$.  As $\si$ ranges over $R_{\pm d}$, these affines
cover $X \mod \pm$, so the result holds globally.  \fp

\begin{s}{Theorem}
\label{2d}
Suppose $X \mod + \neq \varnothing \neq X \mod -$.  Then the blow-ups
of $X \mod \pm$ at $\ci_\pm \mod \pm$, and the blow-up of $X \mod 0$
at $(\ci^+ + \ci^-) \mod 0$, are all naturally isomorphic to the
irreducible component of the fibred product ${X \mod -} \times_{X \mod
0} {X \mod +}$ dominating $X \mod 0$.
\end{s}

\pf.  The blow-up of $X \mod +$ at $\ci_+ \mod +$ is $\Proj
\bigoplus_n H^0((\ci_+ \mod +)^n(dn))$ for $d$ sufficiently divisible.
But by \re{2b} $H^0((\ci_+ \mod +)^n(dn)) = R_{(d+1)n, n}$, so the
blow-up is $\Proj \bigoplus_n R_{(d+1)n, n}$.  There is a surjection
of $R_0$-modules $R_{(d+1)n} \otimes_{R_0} R_{-n} \to R_{(d+1)n, n}$,
so the blow-up embeds in $\Proj \bigoplus_n R_{(d+1)n} \otimes_{R_0}
R_{-n}$.  This is precisely the fibred product ${X \mod -} \times_{X
\mod 0} {X \mod +}$, with polarization $\co(d+1,1)$.  By \re{2m} this
naturally contains $(X \sans X^\pm) \mod \pm$ as a nonempty open set,
but so does the blow-up.  The blow-up is certainly irreducible, so it
is the component containing $(X \sans X^\pm) \mod \pm$.  The proof for
$X \mod -$ is similar.  Likewise, the blow-up of $X \mod 0$ at $(\ci^+
+ \ci^-) \mod 0$ is $\Proj \bigoplus_n (\ci^+ + \ci^- \mod 0)^n =
\Proj \bigoplus_n R_{dn,dn} = \Proj \bigoplus_n R_{n,n}$.  This embeds
in the fibred product with polarization $\co(1,1)$, but by \re{2m}
contains $(X \sans X^\pm) \mod \pm$ as a nonempty open set.  \fp \bl

The ideal sheaves $\ci^\pm$ are supported on $X^\pm$, so $\ci^\pm \mod
\pm$ are supported on $X^\pm \mod \pm$.  But they are not just
$\ci_{X^\pm \mod \pm}$, as the following counterexample shows.

\begin{s}{Counterexample}
\label{2g}
To show that the ideal sheaves of $X^\pm \mod \pm$ cannot
generally replace $\ci_\pm \mod \pm$ in \re{2d}.
\end{s}

In other words, the blow-up may be weighted, not just the usual
blow-up of a smooth subvariety.  Let $G$ act on $X = \Aff^3$ by
$\la(w,x,y) = (\la^{-1} w, \la x, \la^2 y)$; in other words, $w \in
R_{-1}$, $x \in R_1$, and $y \in R_2$.  Then $X \mod 0 = \Spec
\k[wx,w^2y] = \Aff^2$, and
$X \mod - = \Proj \k[wx,w^2y, zw] = \Aff^2$, where the $\N$-grading of
every variable is 0 except $z$, which is graded by 1.  However,
\beqas
X \mod + & = & \Proj \k[wx,w^2y,zx,z^2y] \\
& = & \Proj \k[wx,w^2y,z^2(w^2x^2),z^2(w^2y)] \\
& = & \Proj \k[u,v,zu^2,zv],
\eeqas
which is $\Aff^2$ blown up at the ideal $\langle u^2,v \rangle$.  This
has a rational double point, so is not the usual blow-up at a point.
\fp

The paper \cite{bp} of Brion and Procesi asserts (in section 2.3) a
result very similar to \re{2d}.  They state that the two quotients are
related by blow-ups---``\'eclatements''---of certain subvarieties.
The counterexample  above shows that the blow-ups must sometimes be
weighted, that is, must have non-reduced centres.  Brion and Procesi
do not state this explicitly, but they are no doubt aware of it.
Another minor contradiction to their result is furnished by the
following counterexample.

\begin{s}{Counterexample}
To show that the fibred product of \re{2d} can be reducible, and the
blow-up one of its irreducible components.
\end{s}

Let $G$ act on the singular variety $X = \Spec \k[a^2,ab,b^2,c,d]/
\langle ad-bc \rangle$ where $a,b,c,d$ are of degree 1, acted on with
weights $1,-1,1,-1$ respectively.  Then
\beqas
X \mod 0 & = & \Spec \k[ab,cd,a^2d^2,b^2c^2]/\langle ad-bc \rangle  \\
& = & \Spec \k[ab,cd] \\
& = & \Aff^2.
\eeqas
But, using the same $\N$-grading convention as in the previous example,
\beqas
X \mod - & = & \Proj \k[ab,cd,z^2b^2,zd] \\
& = & \Proj \k[u,v,zu,zv] \\
& = & \blowup,
\eeqas
that is, the blow-up of $\Aff^2$ at the origin, and by symmetry $X
\mod + \cong \blowup$ as well.  Taking $d =2$ gives ideal sheaves
$\ci_+ = \langle a^2,c^2 \rangle$ and $\ci_- = \langle b^2,d^2
\rangle$; the sheaf of invariants of both is$\langle u^2,v^2 \rangle$,
the ideal sheaf of twice the exceptional divisor.  Hence blowing up
$\ci_\pm \mod \pm$ does nothing.  The fibred product, however, is
$\blowup \times_{\Aff^2} \blowup$, which is not just $\blowup$: it has
two components, isomorphic to $\blowup$ and $\Pj^1 \times \Pj^1$
respectively and meeting in a $\Pj^1$.  \fp \bl

Following Bialynicki-Birula \cite{bb}, define the {\em trivial
$w_i$-fibration} over an affine variety $Y$ to be $\Aff^r \times Y$,
with a $G$-action induced by the action on $\Aff^r$ with weights $w_i$.
A {\em $w_i$-fibration} over $Y$ is a variety over $Y$, with a
$G$-action over the trivial action on $Y$, which is locally the
trivial $w_i$-fibration.  As Bialynicki-Birula points out, a
$w_i$-fibration need not be a vector bundle, because the transition
functions need not be linear.  But if all the $w_i$ are equal, then it
is a vector bundle.

Suppose now that $G$ acts on an affine variety $X$ which is {\em
smooth}.  Then it will be proved in (a) below that $X^0$ is also
smooth.  Purely for simplicity, suppose that it is also connected.
(If not, the following theorem is still valid, but the fibrations
involved need have only locally constant rank and weights.)  The group
$G$ acts on the normal bundle $N_{X^0 / X}$.  Let $N^\pm_X$, or simply
$N^\pm$, be the subbundles of positive and negative weight spaces for
this action, with weights $w^\pm_i \in \Z$.

\begin{s}{Theorem}
\label{2e}
Suppose $X$ is smooth.  Then \mbox{\rm (a)} $X^0$ is smooth; \mbox{\rm
(b)} as varieties with $G$-action, $\pi_\pm: X^\pm \to X^0$ are
naturally $w^\pm_i$-fibrations; \mbox{\rm (c)} $N_{X^0 / X}$ has no
zero weights, so equals $N^+ \oplus N^-$; \mbox{\rm (d)} the normal
bundles $N_{X^0 / X^\pm} = N^\pm$; \mbox{\rm (e)} if all $w^\pm_i =
\pm w$ for some $w$, then $\ci^\pm$ and $\ci^+ +\ci^-$ cut out the
$d/w$th infinitesimal neighbourhoods of $X^\pm$ and $X^0$ respectively.
\end{s}

Parts (b), (c), (d) are the {\em Bialynicki-Birula decomposition
theorem} \cite[Thm.\ 4.1]{bb}.  Another version of this result,
possibly more familiar, gives a Morse-style decomposition of a
projective variety into a disjoint union of $w_i$-fibrations.  It
follows easily from this \cite{bb}.

\pf.  First consider the case of a finite-dimensional vector space
$V$, acted on linearly by $G$.  Then $V = \Spec S$ for $S$ a
$\Z$-graded polynomial algebra.  This decomposes naturally into three
polynomial algebras, $S = S^- \otimes S^0 \otimes S^+$, corresponding
to the subspaces of negative, zero, and positive weight.  Then
$\ci^\pm = \langle S^\mp_{\mp d} \rangle$, $V^\pm = \Spec S^\pm
\otimes S^0$, and $V^0 = \Spec S^0$.  Parts (a)--(e) all follow
easily.  Indeed, the fibrations are naturally trivial.

To return to the general case, note first that if $U \subset X$ is a
$G$-invariant open affine, then $U = \Spec F^{-1}R$ for some $F
\subset R_0$.  Hence $(F^{-1}R)_i = F^{-1}(R_i)$ for each $i$, so
$\ci^\pm_U = \ci^\pm_X |_U$, $U^\pm = X^\pm \cap U$, $U^0 = X^0 \cap
U$, and $\pi_\pm$ is compatible with restriction.  Consequently, the
whole theorem is local in the sense that it suffices to prove it for a
collection of $G$-invariant open affines in $X$ containing $X^+ \cup
X^-$.

Now for any closed point $x \in X^0$, apply the Luna slice theorem
\cite{luna,mf} to $X$.  Since $G_x = G$, the Luna slice is a
$G$-invariant open affine $U \subset X$ containing $x$, and $G
\times_{G_x} N_x = T_x X$.  Hence there is a strongly \' etale
$G$-morphism (see \cite{luna,mf}) $\phi: U \to V = T_x X$ such that
$\phi(x) = 0$.  In particular, $U = U \mod 0 \times_{V \mod 0} V$.
Any $y \in X^+ \cup X^-$ is contained in some such $U$: indeed, just
take $x = \pi_\pm (y)$.  It therefore suffices to prove the theorem
for $U$, so by abuse of notation, say $U = \Spec R$ from now on.  The
$G$-morphism $\phi$ then corresponds to a graded homomorphism $S \to
R$, where $S$ is a $\Z$-graded polynomial ring, such that $R = R_0
\otimes_{S_0} S$.  In particular, $R_i = R_0 \otimes_{S_0} S_i$ for
each $i$.  Hence $R_{\pm d}$ and $S_{\pm d}$ generate the same ideals
in $R$, so $\ci^\pm_U = \phi^{-1}\ci^\pm_V$.  Also, $U^\pm = \phi^{-1}
V^\pm = V^\pm \times_{V \mod 0} U \mod 0$, and $U^0 = \phi^{-1}V^0 =
V^0 \times_{V \mod 0} U \mod 0$.  This immediately implies (a).  Since
$V^\pm \to V^0$ are trivial fibrations, it also gives the local
trivialization of $X^\pm \to X^0$ near $x$ needed to prove (b).

Parts (c) and (d) also follow, since $\phi$, as an \' etale
$G$-morphism, satisfies $\phi^* N_{V^0/V} = N_{U^0/U}$ as bundles with
$G$-action, so in particular $\phi^* N^\pm_V = N^\pm_U$.  The
hypotheses of part (e) imply that $d$ is a multiple of $w$; the
conclusion holds if and only if the map $R_{\mp w}^{d/w} \to R_{\mp
d}$ is surjective.  This is true for $S$, and follows in general from
$R_i = R_0 \otimes_{S_0} S_i$.  \fp

The above methods can be used to describe the local structure of the
non-reduced schemes cut out by $\ci^\pm$ even when not all $w^\pm_i =
\pm w$, but we will not pursue this.

\begin{s}{Corollary}
\label{2j}
Suppose $X$ is smooth.  Then $X^\pm \mod \pm$ are locally trivial
fibrations over $X^0$ with fibre the weighted projective space
$\Pj(|w_i^\pm|)$.  \fp
\end{s}

If $X \mod + \neq \varnothing \neq X \mod -$, these are the supports
of the blow-up loci of \re{2d}.  On the other hand, if $X \mod - =
\varnothing$, then $X^+ = X$ and $X^0 = X \mod 0$, so this says the
natural morphism $X \mod + \to X \mod 0$ is a weighted projective
fibration.

\begin{r}{Remark}
\label{2f}
It follows from the above corollary that, in this smooth case, ${X^-
\mod -} \times_{X \mod 0} {X^+ \mod +}$ is irreducible of codimension
1 in ${X \mod -} \times_{X \mod 0} {X \mod +}$.  It must therefore be
exactly the exceptional divisor of each of the two blow-ups of
\re{2d}.  In other words, the latter fibred product is irreducible,
and is isomorphic to each of the blow-ups.  This implies that, when
$X$ is smooth, the surjection of $R_0$-modules $R_{i} \otimes_{R_0}
R_{-j} \to R_{i,j}$ is an isomorphism for $i,j > 0$ sufficiently
divisible.  However, I do not know of a direct algebraic proof of this
fact.
\end{r}

\begin{s}{Theorem}
\label{2k}
Suppose $X$ is smooth, and that all $w_i^\pm = \pm w$ for some $w$.
Then $X^\pm \mod \pm$ are naturally isomorphic to the projective
bundles $\Pj(N^\pm)$ over the fixed-point set $X^0$, their normal
bundles are naturally isomorphic to $\pi_\pm^* N^\mp(-1)$, and the
blow-ups of $X \mod \pm$ at $X^\pm \mod \pm$, and of $X \mod 0$ at
$X^0 \mod 0$, are all naturally isomorphic to the fibred product ${X
\mod -} \times_{X \mod 0} {X \mod +}$.
\end{s}

Note that if each 0-dimensional stabilizer on $X$ is trivial, then all
$w^\pm_i = \pm 1$.

\pf.  All the blow-ups and the fibred product are empty if either $X
\mod +$ or $X \mod -$ is empty, so suppose they are not.  By the
observation of Bialynicki-Birula quoted above, if all $w^\pm_i = \pm
w$, then the fibrations $X^\pm \to X^0$ are actually vector bundles.
But any vector bundle is naturally isomorphic to the normal bundle of
its zero section, so by \re{2e}(d) $X^\pm \mod \pm = \Pj(X^\pm) =
\Pj(N^\pm)$, and the natural $\co(1)$ bundles correspond.  By
\re{2e}(e), $\ci^\pm$ cut out the $d/w$th infinitesimal neighbourhoods
of $X^\pm$.  This means that $R_{\mp w}^{d/w} \to R_{\mp d}$ are
surjective and hence that $\ci^\pm \mod \pm$ and $(\ci^+ + \ci^-) \mod
0$ cut out the $d/w$th infinitesimal neighbourhoods of $X^\pm \mod
\pm$ and $X^0 \mod 0$ respectively.  Since blowing up a subvariety has
the same result as blowing up any of its infinitesimal neighbourhoods,
the result follows from \re{2d} and \re{2f}, except for the statement
about normal bundles.  To prove this, recall first that if $E$ is the
exceptional divisor of the blow-up $\tilde{Y}$ of any affine variety
$Y$ at $Z$, then $N_{E/\tilde{Y}} = \co_E(-1)$, and $N^*_{Z/Y} =
(R^0\pi) N^*_{E/\tilde{Y}}$.  Applying this to the case $Y = X \mod 0$
shows that the normal bundle to ${X^- \mod -} \times_{X^0} {X^+ \mod
+}$ is the restriction of $\co(-1,-1) \to {X \mod -} \times_{X \mod 0}
{X \mod +}$, which is exactly the obvious $\co(-1,-1) \to \Pj N^+
\times_{X^0} \Pj N^-$.  The normal bundle to $\Pj N^\pm$ is then just
$(R^0\pi_\mp \co(1,1))^* = \pi_\pm^* N^\mp(-1)$.  \fp

\begin{r}{Example}
The simplest non-trivial example of these phenomena is also the
best-known; indeed it goes back to an early paper of Atiyah \cite{a}.
Let $X = \Aff^4$, and let $G$ act by $\lambda \cdot(v,w,x,y)=
(\lambda v, \lambda w, \lambda^{-1} x, \lambda^{-1} y)$.  Then $X \mod
0 = \Spec \k[vx,vy,wx,wy] = \Spec[a,b,c,d]/\langle ad-bc \rangle$, the
affine cone on a smooth quadric surface in $\Pj^3$.  But, using the
$\N$-grading conventions of the previous examples,
\beqas
X \mod + & = & \Proj \k[vx,vy,wx,wy,zv,zw] \\
 & = & \Proj \k[a,b,c,d,za,zc] / \langle ad-bc \rangle.
\eeqas
This is the blow-up of $X \mod 0$ at the Weil divisor cut out by $a$
and $c$.  But this Weil divisor is generically Cartier, so the
blow-down $X \mod + \to X \mod 0$ is generically an isomorphism even
over the divisor.  The exceptional set of the morphism therefore has
codimension 2; indeed, it is the $\Pj^1$ lying over the cone point.
Likewise, $X \mod - = \Proj \k[a,b,c,d,zb,zd]/ \langle ad-bc \rangle$,
and similar remarks apply by symmetry.  By \re{2k} the fibred product
${X \mod -} \times_{X \mod 0} {X \mod +}$ is the common blow-up of $X
\mod \pm$ at these $\Pj^1$, and also the blow-up of $X \mod 0$ at the
cone point.  This is exactly the proper transform of the quadric cone
in $\Aff^4$ blown up at the origin, so it has fibre $\Pj^1 \times
\Pj^1$ over the cone point, as expected.
\end{r}


\bit{The space of linearizations}

In \S\S3, 4 and 5 we will generalize the results of \S1 in three
directions.  First, the group $G$ may now be any reductive group over
$\k$.  Second, the variety $X$ may now be any quasi-projective variety
over $\k$, projective over an affine variety.  Finally, the
linearization may be arbitrary.  Before doing this, though, we will
prove some general facts, in the case where $X$ is normal and
projective, about the structure of the group of all linearizations.
This will show how to apply our general results in this case.

So in this section, suppose $X$ is normal and projective over $\k$,
and that $G$ is a reductive group over $\k$ acting on $X$.  We first
recall a few familiar facts about the Picard group.

In the Picard group $\Pic$ of isomorphism classes of line bundles, the
property of ampleness depends only on the algebraic equivalence class
of the bundle.  Hence there is a well-defined ample subset $\A$ of the
N\'eron-Severi group $\NS$ of algebraic equivalence classes of line
bundles.  This determines the {\em ample cone} $\A_\Q = \A
\otimes_{\N} \Q_{\geq 0} \subset \NS_\Q = \NS \otimes \Q$.  The
N\'eron-Severi group is finitely-generated, so $\NS_\Q$ is a
finite-dimensional rational vector space.  We will refer to an element
of $\A$ as a {\em polarization}, and an element of $\A_\Q$ as a {\em
fractional polarization}.

Now let $\Pic^G$ be the group of isomorphism classes of linearizations
of the $G$-action (cf.\ 1, \S3 of GIT).  There is a forgetful
homomorphism $f: \Pic^G \to \Pic$, whose kernel is the group of
linearizations on $\co$, which is exactly the group $\chi(G)$ of
characters of $G$.  There is an equivalence relation on $\Pic^G$
analogous to algebraic equivalence on $\Pic$; it is defined as
follows.  Two linearizations $L_1$ and $L_2$ are said to be {\em
$G$-algebraically equivalent} if there is a connected variety $T$
containing points $t_1, t_2$, and a linearization $L$ of the
$G$-action on $T \times X$ induced from the second factor, such that
$L|_{t_1 \times X} \cong L_1$ and $L|_{t_2 \times X} \cong
L_2$.

\begin{s}{Proposition}
\label{1d}
If $L$ is an ample linearization, then $X^{ss}(L)$, and the quotient
$X \mod G(L)$ regarded as a polarized variety, depend only on the
$G$-algebraic equivalence class of $L$.
\end{s}

\pf. The statement about $X^{ss}(L)$ is proved exactly like Cor.\ 1.20
of GIT, except that the Picard group $P$ is replaced by $T$.  The
statement about $X \mod G(L)$ as a variety then follows from this,
since $X \mod G(L)$ is a good quotient of $X^{ss}(L)$, hence a
categorical quotient of $X^{ss}(L)$, so is uniquely determined by
$X^{ss}(L)$.  As for the polarization, note that, if $L_1$ and $L_2$
are $G$-linearly equivalent ample linearizations, then the
linearization $L$ on $T \times X$ inducing the equivalence can be
assumed ample: just tensor $L$ by the pullback of a sufficiently ample
bundle on $T$.  Then $\co(1) \to (T \times X) \mod G (L)$ is a
family of line bundles on $X \mod G$ including $\co(1) \to X \mod G
(L_1)$ and $\co(1) \to X \mod G (L_2)$, so these are algebraically
equivalent.  \fp

So define $\NS^G$ to be the group of $G$-algebraic equivalence classes
of linearizations.  In light of \re{1d}, by abuse of terminology an
element of $\NS^G$ will frequently be called just a linearization.
The forgetful map $f$ descends to $f: \NS^G \to \NS$.

\begin{s}{Proposition}
\label{2q}
This new $f$ has kernel $\chi(G)$ modulo a torsion subgroup.
\end{s}

\pf.  Let $M \to \Pic_0 X \times X$ be the Poincar\'e line bundle, and
let $G$ act on $\Pic_0 X \times X$, trivially on the first factor.  By
Cor.\ 1.6 of GIT, some power $M^n$ of $M$ admits a linearization.
Since the $n$th power morphism $\Pic_0 X \to \Pic_0 X$ is surjective,
this shows that any element of $\Pic_0 X$ has a linearization
$G$-algebraically equivalent to a linearization on $\co$.  Hence $\ker
f$ is $\chi(G)$ modulo the subgroup of linearizations on $\co$
which are $G$-algebraically equivalent to the trivial linearization.
We will show that this subgroup is torsion.

Suppose there is a linearization $L_1$ on $\co$ which is
$G$-algebraically equivalent to the trivial linearization.  Then there
exist $T$ containing $t_1$, $t_2$ and $L$ as in the definition of
$G$-algebraic equivalence.  There is an induced morphism $g: T \to
\Pic_0 X$ such that $t_1, t_2 \mapsto \co$.  As before, let $M^n$ be
the power of the Poincar\'e bundle admitting a linearization.  Then $N
= (1 \times g)^* M^n \otimes L^{-n}$ is a family of linearizations on
$\co \to X$.  Since the isomorphism classes of such linearizations
form the discrete group $\chi(G)$, $L_1^n = N_{t_1}^{-1} \otimes
N_{t_2}$ is trivial as a linearization.  \fp \bl

Hence $\NS^G$ is finitely-generated and $\NS^G_\Q = \NS^G \otimes
\Q$ is again a finite-dimensional rational vector space.  We refer to
an element of $\NS^G_\Q$ as a {\em fractional linearization}.

The map $f: \NS^G \to \NS$ is not necessarily surjective (see 1, \S3
of GIT).  But by Cor.\ 1.6 of GIT, $f_\Q: \NS^G_\Q \to \NS_\Q$ is
surjective.  By \re{2q}, the kernel is $\chi(G) \otimes \Q$, the group
of {\em fractional characters}.  (Not to be confused with $f_\Q$ is
the natural surjective linear map $\NS^G_\Q(X) \to \NS_\Q(X \mod G)$:
this is induced by descent, since divisor classes always descend over
$\Q$.)

An ample linearization $L$ is said to be {\em $G$-effective}
if $L^n$ has a $G$-invariant section for some $n > 0$.  This is
equivalent to having a semistable point, so \re{1d} shows that
$G$-effectiveness depends only on the $G$-algebraic equivalence class
of the linearization.  Hence there is a well-defined $G$-effective
subset $\E^G \subset f^{-1}\A \subset \NS^G$, and a {\em $G$-effective
cone} $\E^G_\Q = \E^G \otimes_{\N} \Q_{\geq 0} \subset \NS^G_\Q$.

Now a linearization $L$ determines a quotient $X \mod G$ if $L$ is
ample; then $X \mod G \neq \varnothing$ if and only if $L$ is also
$G$-effective.  Of course, we can also tensor by $\Q$, allowing
fractional linearizations; the quotient $X \mod G$ will then be
fractionally polarized.  Hence any fractional linearization $L \in
f_\Q^{-1}(\A_\Q) \subset \NS^G_\Q$ defines a fractionally
polarized quotient, which will be nonempty if and only if $L \in
\E^G_\Q$.  Replacing a fractional linearization $L$ by $L^n$ for
some positive $n \in \Q$ has no effect on the quotient, except to
replace the fractional polarization $\co(1)$ by $\co(n)$. \bl

The first result describing the dependence of the quotient $X \mod
G(L)$ on the choice of $L \in \NS^G_\Q$ is the following, which is
analogous to the Duistermaat-Heckman theory in symplectic geometry.

\begin{s}{Theorem}
\label{1b}
The $G$-effective cone $\E_\Q^G$ is locally polyhedral in the ample
cone $f_\Q^{-1} \A_\Q$.  It is divided by homogeneous {\rm walls},
locally polyhedral of codimension 1 in $f_\Q^{-1} \A_\Q$, into convex
{\rm chambers} such that, as $t$ varies within a fixed chamber, the
semistable set $X^{ss}(t)$, and the quotient $X \mod G(t)$, remain
fixed, but $\co(1)$ depends affinely on $t$. If $t_0$ is on a wall or
walls, or on the boundary of $\E_\Q^G$, and $t_+$ is in an adjacent
chamber, then there is an inclusion $X^{ss}(+) \subset X^{ss}(0)$
inducing a canonical projective morphism $X \mod G(+) \to X \mod
G(0)$.
\end{s}

The proposition above could be proved directly, using Kempf's descent
lemma \cite{dn} for the statement about $\co(1)$, and Mumford's
numerical criterion \cite[4.9; GIT Thm.\ 2.1]{n} for the rest.  But it
will follow easily from the construction \re{1h} to be introduced in
the next section, so we put off the proof until then. \bl

Theorem \re{1b} asserts that the walls are locally polyhedral, and in
particular, locally finite; but with a little more effort we can prove
a global result.

\begin{s}{Theorem}
\label{1s}
There are only finitely many walls.
\end{s}

\pf.  Suppose not.  Then there exists an infinite sequence $\{ C_i \}$
of chambers such that, for any $m,n \geq 0$, the convex hull of $C_n
\cup C_{n+m}$ intersects the interior of $C_{n+1}$ nontrivially.
Indeed, let $C_0$ be any chamber; then there exists a wall $W_0$
bounding it such that there are infinitely many chambers on the other
side of $W_0$ (or more properly, the affine hyperplane containing
$W_0$).  Let $C_1$ be the other chamber bounded by $W_0$.
Inductively, given $C_0, \dots, C_n$ such that $C_n$ is on the other
side of $W_i$ from $C_i$ for all $i < n$, there is a wall $W_n$
bounding $C_n$ such that there are infinitely many chambers which for
all $i \leq n$ are on the other side of $W_i$ from $C_i$.  Let
$C_{n+1}$ be the other chamber bounded by $W_n$.  For a sequence
chosen in this manner, $C_{n+m}$ is on the other side of $W_n$ from
$C_n$, so the convex hull of $C_n \cup C_{n+m}$ meets the interior of
$C_{n+1}$.

Choose an $L_i$ in the interior of each $C_i$.  For any fixed $x \in
X$, the set $\{ L \in \NS_\Q \st x \in X^{ss}(L) \}$ is convex, since
$s_\pm \in H^0(L_\pm)^G$, $s_\pm(x) \neq 0$ imply $s_+ \cdot s_- \in
H^0(L_+ \otimes L_-)^G$, $(s_+ \cdot s_-)(x) \neq 0$.  But by \re{1b}
it is also a union of chambers, so by induction it includes $C_n \cup
C_{n+m}$ if and only if it includes $C_{n+i}$ for all $i \leq m$.  Its
intersection with $\{ L_i \}$ is therefore the image of an interval in
$\N$.  Hence $X^{ss}(L_{i+1}) \sans X^{ss}(L_i)$ are all disjoint; but
each one is open in the complement of $X^{ss}(L_0)$ and the preceding
ones.  Since varieties are noetherian, this implies there exists $i_0$
such that for all $i \geq i_0$, $X^{ss}(L_{i+1}) \sans X^{ss}(L_i) =
\varnothing$, and hence $X^{ss}(L_{i+1}) \subset X^{ss}(L_i)$.  There
is therefore an infinite sequence of dominant projective morphisms $$
\cdots \to X \mod G(L_{i_0 + 2}) \to X \mod G(L_{i_0 + 1}) \to X\mod
G(L_{i_0}).$$ Hence the N\'eron-Severi group of $X \mod G(L_i)$ has
arbitrarily large rank for some $i$.  But as mentioned before, there
is a natural surjective linear map $\NS^G_\Q(X) \to \NS_\Q(X \mod
G(L_i))$ for all $i$.  Since $\NS^G_\Q(X)$ is finite-dimensional, this
is a contradiction.  \fp

\bit{The general case}

We now embark on our generalization of the results of \S1.  So let $G$
be a reductive group over $\k$, acting on a quasi-projective variety
$X$ over $\k$, projective over an affine variety.  This is the largest
category in which geometric invariant theory guarantees that the
semistable set has a good quotient.

All of the arguments in this section use the following trick.

\begin{r}{Construction}
\label{1h}
Let $L_1 , \dots , L_{n+1}$ be ample linearizations.  Let $\Delta$ be
the $n$-simplex $\{ (t_i) \in \Q^{n+1} \st \sum t_i = 1 \} $.  Then
for any $t = (t_i) \in \Delta$, $L(t) = \bigotimes_i L_i^{t_i}$ is a
fractional linearization on $X$.  We refer to the set $\{ L(t) \st t
\in \Delta \}$ as an {\em $n$-simplicial family} of fractional
linearizations.

Put
$$Y = \Pj(\bigoplus_i L_i) = \Proj \sum_{j_i \in \N} H^0(\bigotimes_i
L_i^{j_i}),$$
and let $q: Y \to X$ be the projection.  Then $G$ acts naturally on
$\bigoplus_i L_i$, hence on $Y$ with a canonical linearization on the
ample bundle $\co(1)$.  Likewise, the $n$-parameter torus $T = \{\la
\in \k^{n+1} \st \prod_i \la_i = 1 \}$ acts on $\bigoplus_i L_i$ by
$\la(u_i) = (\la_i u_i)$, and hence on $Y$.  This $T$-action commutes
with the $G$-action.  Moreover, since it comes from $\bigoplus_i L_i$,
it too is linearized on $\co(1)$.  But this obvious linearization is
not the only one.  Indeed, any $t \in \Delta$ determines a fractional
character of $T$ by $t(\la) = \prod_i \la_i^{t_i}$; then $\la(u_i) =
(t^{-1}(\la) \la_i u_i)$ determines a fractional linearization
depending on $t$.  This gives an $n$-simplicial family $M(t)$ of
fractional linearizations on $\co(1)$ of the $T$-action on $Y$, each
compatible with the canonical linearization of the $G$-action.  In
other words, $M(t)$ is a family of fractional linearizations of the $G
\times T$-action on $Y$.  Let $Y^{ss}(t)$ be the semistable set for
this action and linearization, and let $Y^{ss}(G)$ be the semistable
set for the $G$-action alone.  For any $t$, $Y^{ss}(t) \subset
Y^{ss}(G)$.

With respect to $M(t)$, $T$ acts trivially on $H^0(\bigotimes_i
L_i^{j_i})$ if and only if $j_i = m t_i$ for some fixed $m$.  Hence
the subalgebra of $T$-invariants is $\sum_m H^0 ((\bigotimes_i
L_i^{t_i})^m)$.  The quotient $Y \mod T(t)$ is therefore $X$, but with
the residual $G$-action linearized by $L(t)$.  Consequently $(Y \mod
T(t)) \mod G = X \mod G(t)$, the original quotients of interest.
Moreover, $X^{ss}(t) = q(Y^{ss}(t))$.

Let $Z$ be the quotient $Y \mod G$ with respect to the canonical
linearization defined above, and let $p: Y^{ss}(G) \to Z$ be the
quotient.  Then the $M(t)$ descend to an $n$-simplicial family $N(t)$
of fractional linearizations on $\co(1)$ of the residual $T$-action on
$Z$, and $Y^{ss}(t) = q^{-1}(Z^{ss}(t))$.  When two group actions
commute, the order of taking the quotient does not matter, so $(Y \mod
T(t)) \mod G = (Y \mod G) \mod T(t) = Z \mod T(t)$.  So we have
constructed a variety $Z$, acted on by a torus $T$, and a simplicial
family $N(t)$ of fractional linearizations on $\co(1)$, such that $Z
\mod T(t) = X \mod G(t)$.  Moreover, $X^{ss}(t) = q(p^{-1}(Z^{ss}(t)))$.
\end{r}

As a first application of this construction, let us prove the result
asserted in the last section.

\pf\ of \re{1b}.  The result is relatively easy in the case where $X =
\Pj^n$ and $G$ is a torus $T$.  Indeed, the $T$-effective cone is
globally polyhedral, as is each chamber; for details see \cite{bp,
toric}.  In the general case, choose a locally finite collection of
simplices in $f_\Q^{-1}\A_\Q \subset \NS^G_\Q$ such that every vertex
is in $\NS^G$, and for every $L \in f_\Q^{-1}\A_\Q$, some fractional
power $L^m$ is in one of the simplices.  By the homogeneity property
mentioned just before the statement of \re{1b}, it suffices to prove
the statement analogous to \re{1b} where $f_\Q^{-1}\A_\Q$ is replaced
by the simplex parametrizing each of these families.  The construction
\re{1h} applies, so there exists $Z \subset \Pj^n$ and a simplicial
family $N(t)$ in $\NS^T(\Pj^n)$ such that $X \mod G(t) = Z \mod T(t)$
for all $t \in \Delta$.  The conclusions of the theorem are preserved
by restriction to a $T$-invariant subvariety, so they hold for $Z$ and
$Z \mod T(t)$, and hence for $X$ and $X \mod G(t)$.  \fp \bl

The rest of this section and all of \S\S4 and 5 will refer to the
following set-up.  Let $X$ and $G$ be as before.  Let $L_+$ and $L_-$
be ample linearizations such that, if $L(t) = L_+^t L_-^{1-t}$ for $t
\in [-1,1]$, there exists $t_0 \in (-1,1)$ such that $X^{ss}(t) =
X^{ss}(+)$ for $t > t_0$ and $X^{ss}(t) = X^{ss}(-)$ for $t < t_0$.

For example, \re{1b} implies that this is the case if $X$ is normal
and projective, $L_\pm$ are in adjacent chambers, and the line segment
between them crosses a wall only at $L(t_0)$.  Even in the normal
projective case, however, there are other possibilities; for example,
$L_\pm$ could both lie in the same wall, or $L(t_0)$ could lie on the
boundary of $\E^G_\Q$.  In future, $L(t_0)$ will be denoted $L_0$. \bl

The following lemma shows how to globalize the results of \S1 within
this set-up.  Suppose $T \cong \kst$ acts on $X$, and let $\si
\in H^0(X, L_0^n)^T$ for some $n$, so that $X_\si \subset
X^{ss}(0)$ is a $T$-invariant affine.

\begin{s}{Lemma}
\label{1e}
Suppose $f(L_-) = f(L_+)$.  Then \mbox{\rm (a)} $X^{ss}(\pm) \subset
X^{ss}(0)$; \mbox{\rm (b)} $X_\si^{ss}(0) = X_\si \cap X^{ss}(0)$; and
\mbox{\rm (c)} $X_\si^{ss}(\pm) = X_\si \cap X^{ss}(\pm)$; so there is
a natural commutative diagram
$$\begin{array}{ccc}
     X_\si \mod \pm & \emb  & X \mod \pm \vspace{.7ex} \\
     \down{} &  & \down{} \vspace{.7ex} \\
     X_\si \mod 0 & \emb  & X \mod 0
\end{array} $$
whose rows are embeddings.
\end{s}

\pf.  Put $R_m = H^0(X, L_0^m)$, so that $X = \Proj \bigoplus_{m \in
\N}R_m$, and let $R_m = \bigoplus_{n \in \Z} R_{m,n}$ be the weight
decomposition for the $\kst$-action.  Suppose $x \in X^{ss}(+) \sans
X^{ss}(0)$.  Then for $m>0$, every element of $R_{m,0}$ vanishes at
$x$.  Since $\bigoplus_m R_m$ is finitely-generated, this implies
that, for $m/n$ large enough, every element of $R_{m,n}$ vanishes at
$x$.  But then there exists $t>t_0$ such that $x \notin X^{ss}(t)$,
contradicting the set-up.  The proof for $X^{ss}(-)$ is similar.  This
proves (a).

Without loss of generality suppose $\si \in H^0(X, L_0)$.  Then $\si
\in R_{1,0}$ and $X_\si = \Spec (\si^{-1}R)_0$.  Since $\bigoplus_m
R_{m,0}$ is finitely-generated, for $m$ large the map $R_{m,0} \to
(\si^{-1} R)_{0,0}$ given by dividing by $\si^m$ is surjective.  But
$R_{m,0} = H^0(X, L_0^m)^T$ and $(\si^{-1} R)_{0,0} = H^0(X_\si,
L_0^m)^T$ (the latter since $L_0^m$ is trivial on $X_\si$), so this
implies that $X_\si^{ss}(0) = X_\si \cap X^{ss}(0)$, hence that $X_\si
\mod 0$ embeds in $X \mod 0$.  This proves (b).

Without loss of generality take $L_+$ to be $L_0$ twisted by the
fractional character $\la \mapsto \la^{1/p}$ for $p$ large.  Since $R$ is
finitely-generated, $R_{m,n} \to (\si^{-1} R)_{0,n}$ is
surjective for $m/n$ large.  But for $m = np$, $R_{m,n} = H^0(X,
L_+^m)^T$ and $(\si^{-1} R)_{0,n} = H^0(X_\si, L_+^m)^T$,
so this implies that $X_\si^{ss}(+) = X_\si \cap X^{ss}(+)$,
hence that $X_\si \mod +$ embeds in  $X \mod +$.  The case
of $L_-$ is similar.  This proves (c).  \fp

Hence, in studying $\kst$-quotients where $L_+ \cong L_-$ as bundles,
we may work locally, using the methods of \S1.

\begin{s}{Theorem}
\label{1f}
If $X \mod G(+)$ and $X \mod G(-)$ are both nonempty, then the
morphisms $X \mod G(\pm) \to X \mod G(0)$ are proper and birational.
If they are both small, then the rational map  $X \mod G(-) \dasharrow X
\mod G(+)$ is a flip with respect to $\co(1) \to X \mod G(+)$.
\end{s}

Again, this could be proved directly, by first examining the stable
sets to show birationality, then applying Kempf's descent lemma to the
linearization $L_+$.  But again, we will use the trick.

\pf\ of \re{1f}.  Perform the construction of \re{1h} on $L_+$ and
$L_-$.  This gives a variety $Z$ with an action of $T \cong \kst$ and
a family $N(t)$ of fractional linearizations with $f_\Q(N(t))$ constant
such that $Z \mod T(t) = X \mod G (t)$.  The whole statement is local
over $X \mod G(0)$, so by \re{1e} it suffices to prove it for affines
of the form $Z_\si$, with the $T$-action and fractional
linearizations $N(t)$.  But this is the case considered in \S1, so
\re{2h} and \re{2i} complete the proof.  \fp

There is a converse to \re{1f} analogous to \re{2o}, which we leave to
the reader.

\begin{r}{Application}
For an application, suppose that $X$ is normal and projective.  Choose
any nonzero $M \in \NS_\Q X$, let $L(t) = L \otimes M^t$ and consider
the ray $\{ L(t) \st t \geq 0 \} \subset \NS^G_\Q(X)$.  By \re{1b},
the quotient $X \mod G(t)$ is empty except for $t$ in some bounded
interval $[0,\omega]$, and this interval is partitioned into finitely
many subintervals in whose interior $X \mod G(t)$ is fixed.  But when
a critical value $t_0$ separating two intervals is crossed, there are
morphisms $X \mod G(t_\pm) \to X \mod G(t_0)$, which by \re{1f} are
birational except possibly at the last critical value $\omega$.  Since
the fractional polarization on $X \mod G(t)$ is the image of $L
\otimes M^{-t}$ in the natural descent map $\NS^G_\Q(X) \to \NS_\Q(X
\mod G(t))$, the descents of $M$ to $\Q$-Cartier divisor classes on $X
\mod G(\pm)$ are relatively ample for each morphism $X \mod G(t_+) \to
X \mod G(t_0)$, and relatively negative for each morphism $X \mod
G(t_-) \to X \mod G(t_0)$.  So suppose that each $X \mod G(t_+) \to X
\mod G(t_0)$ is small when $X \mod G(t_-) \to X \mod G(t_0)$ is small,
and that each $X \mod G(t_+) \to X \mod G(t_0)$ is an isomorphism when
$X \mod G(t_-) \to X \mod G(t_0)$ is divisorial.  It then follows that
the finite sequence of quotients $X \mod G(t)$ runs the $M$-minimal
model programme \cite[(2.26)]{ko} on $X \mod G(L)$, where by abuse of
notation $M$ denotes its image in the descent map.
\end{r}

For some $d > 0$, the ideal sheaves $\langle H^0(X, L_\pm^{nd})^G
\rangle$ and $\langle H^0(X, L_\pm^{d})^G\rangle^n$ on $X$ are equal
for all $n \in \N$.  For such a $d$, let $\ci^\pm = \langle H^0(X,
L_\mp^{d})^G \rangle$ (note the reversal of sign), and let $\ci^\pm
\mod G(\pm)$ be the corresponding sheaves of invariants on $X \mod G
(\pm)$.  Also let $(\ci^+ + \ci^-) \mod G(0)$ be the sheaf of
invariants of the ideal sheaf $\ci^+ + \ci^-$ on $X \mod G(0)$.

\begin{s}{Theorem}
\label{1g}
Suppose $X \mod G(+)$ and $X \mod G(-)$ are both nonempty.  Then the
pullbacks of $(\ci^+ + \ci^-) \mod G(0)$ by the morphisms $X \mod
G(\pm) \to X \mod G(0)$ are exactly $\ci^\pm \mod G(\pm)$, and the
blow-ups of $X \mod G(\pm)$ at $\ci^\pm \mod G(\pm)$, and of $X \mod
G(0)$ at $(\ci^+ + \ci^-) \mod G(0)$, are all naturally isomorphic to
the irreducible component of the fibred product $X \mod G(-) \times_{X
\mod G(0)} X \mod G(+)$ dominating $X \mod G(0)$.
\end{s}

\pf.  Construct a variety $Z$ as in the proof of \re{1f}.
Notice that for $d$ large, since $\ci^\pm_X = \langle H^0(X,
L_\pm^d)^G \rangle$ on $X$ and $\ci^\pm_Z = \langle H^0(Z, N_\pm^d)^T
\rangle$ on $Z$, the pullbacks of both $\ci_X^\pm$ and $\ci_Z^\pm$ to
$Y$ are $\ci_Y^\pm = \langle H^0(Y, M_\pm^d)^{G \times T} \rangle$.
Hence $\ci_X^\pm$ and $\ci_Z^\pm$ have the same sheaves of invariants
on the quotients $Z \mod T(t) = X \mod G (t)$.  It therefore suffices
to prove the statement for $Z$ and $N_\pm$.  All statements are local
over $Z \mod T(0)$, so by \re{1e} it suffices to prove them for
affines of the form $Z_\si$.  But this is the case considered in
\S1, so \re{2b} and \re{2d} complete the proof.  \fp \bl

\bit{The smooth case: strong results}

In the next two sections we seek to generalize the other two main
results of \S1, \re{2j} and \re{2k}.  Indeed, we will give two
different generalizations of each.  The generalizations in \S4 make
fairly strong hypotheses, and prove that, as in \S1, $X^\pm \mod
G(\pm)$ are locally trivial over $X^0 \mod G(0)$.  Moreover, the
proofs are quite easy using the tools already at hand.  Those in \S5
relax the hypotheses somewhat, but conclude only that $X^\pm \mod
G(\pm)$ are locally trivial in the \'etale topology.  The proofs
therefore require \'etale covers and are more difficult; in fact we
confine ourselves to a sketch of the \'etale generalization of
\re{2k}.

Let $X$, $G$, and $L_\pm$ be as in \S3.  As in \re{1h}, let $Y =
\Pj(L_+ \oplus L_-)$, let $T$ be the torus acting on $Y$, let $Z = Y
\mod G$, and let $p: Y^{ss}(G) \to Z$ be the quotient morphism.  Fix
the isomorphism $T \cong \kst$ given by projection on the first
factor.  Define $Y^\pm$, $Y^0$, $Z^\pm$, and $Z^0$ similarly to
$X^\pm$ and $X^0$.  Also let $i_\pm : X \to Y$ be the embeddings given
by the sections at 0 and $\infty$.  Write $q : Y \to X$ for the
projection as before, but let $\pi$ denote the restriction of $q$ to
$Y \sans (i_+(X) \cup i_-(X))$. So in particular $\pi^{-1}(X)$ denotes
$Y \sans (i_+(X) \cup i_-(X))$ itself.

\begin{s}{Lemma}
\label{2r}
$X^{ss}(\pm) \subset X^{ss}(0)$.
\end{s}

\pf.  This is true for $Z$ by \re{1e}(a), but $X^{ss}(\pm) =
q(p^{-1}(Z^{ss}(\pm)))$ and $X^{ss}(0) = q(p^{-1}(Z^{ss}(0)))$.  \fp

\begin{s}{Lemma}
\label{1j}
\mbox{\rm (a)} $i_\pm(X) \cap Y^{ss}(0) = \varnothing$; \mbox{\rm (b)}
$i_\pm(X) \cap Y^{ss}(G) = i_\pm(X^{ss}(\pm))$; \mbox{\rm (c)}
$\pi^{-1}(X) \cap Y^{ss}(0) = \pi^{-1}(X) \cap Y^{ss}(G) =
\pi^{-1}(X^{ss}(0))$.
\end{s}

\pf.  For $i_\pm(x)$ to be in $Y^{ss}(0)$, it must certainly be
semistable for the $T$-action on the fibre $q^{-1}(x) = \Pj^1$.  But
in the fractional linearization $M_0$, $T$ acts with nontrivial weight
on both homogeneous coordinates of $\Pj^1$, so any invariant section
of $\co(n)$ for $n > 0$ must vanish both at $0$ and $\infty$.  Hence
$i_\pm(x)$ are unstable, which proves (a).  However, for $i_\pm(x)$ to
be in $Y^{ss}(G)$ requires only that the section of $\co(n)$ which is
nonzero at $x$ be $G$-invariant.  Pushing down by $q$ shows that
$$H^0(Y, \co(n)) = H^0(X, \bigoplus_{j=0}^n L_+^j \otimes L_-^{n-j}) =
\bigoplus_{j=0}^n H^0(X,  L_+^j \otimes L_-^{n-j}),$$
and a section of $\co(n)$ is nonzero at $i_\pm(x)$ if and only if its
projection on $H^0(X, L_\pm^n)$ is nonzero at $x$.  Hence $i_\pm(x)
\in X^{ss}(G)$ if and only if $x \in X^{ss}(\pm)$, which proves (b).

With respect to the fractional linearization $M(t)$, the $T$-invariant
subspace in the above decomposition consists of that $H^0(X, L_+^j
\otimes L_-^{n-j})$ such that $L_+^j \otimes L_-^{n-j}$ is a power of $L(t)$,
and an invariant section is nonzero on $\pi^{-1}(x)$ if and only if
the corresponding element of $H^0(X, L_+^j \otimes L_-^{n-j})$ is nonzero at
$x$.
Hence there is a $G$-invariant section of some $\co(n)$ non-vanishing
on $\pi^{-1}(x)$ if and only if there is a $G$-invariant section of
some $L(t)^n$ non-vanishing at $x$; this implies $\pi^{-1}(x) \subset
Y^{ss}(G)$ if and only if $x \in \cup_t X^{ss}(t)$, which equals
$X^{ss}(0)$ by \re{2r}.  On the other hand, $x \in X^{ss}(0)$ if and
only if there is a $G$-invariant section of some $L_0^n$ non-vanishing
at $x$, and hence a $G \times T$-invariant section of some $M_0^n$
non-vanishing on $\pi^{-1}(x)$, that is, $\pi^{-1}(x) \subset
Y^{ss}(0)$.  This proves (c).  \fp \bl

Let $X^\pm$ and $X^0$ be the intersections with $X^{ss}(0)$ of the
supports of the sheaves $\ci^\pm$ and $\ci^+ + \ci^-$, defined as in
\S3.  Note that this generalizes the definitions of \S1.  Indeed,
\beqas
X^\pm & = & X^{ss}(0) \sans X^{ss}(\mp); \\
X^0   & = & X^{ss}(0) \sans (X^{ss}(+) \cup X^{ss}(-)).
\eeqas

\begin{s}{Lemma}
\label{1u}
\mbox{\rm (a)} $\pi^{-1}X^\pm = Y^\pm = p^{-1}Z^\pm$; \mbox{\rm (b)}
$\pi^{-1}X^0 = Y^0 = p^{-1}Z^0$.
\end{s}

\pf.  These follow immediately from $\pi^* \ci^\pm_X = \ci^\pm_Y = p^*
\ci^\pm_Z$.  \fp \bl

Choose $x \in X^0$; throughout this section, we will assume the
following.

\begin{r}{Hypothesis}
\label{1y}
Suppose that $X$ is smooth at $x \in X^0$, that $G \cdot x$ is closed
in $X^{ss}(0)$, and that $G_x \cong \kst$.
\end{r}

Note that if $G_x \cong \kst$ for {\em all} $x \in X^0$, then an orbit
in $X^0$ cannot specialize in $X^0$, so it is closed in $X^0$ and
hence in $X^{ss}(0)$.  So the second part of the hypothesis is
redundant in this case.  The third part is necessary, as the
counterexample \re{1r} will show.  But it is always true when $G$ is a
torus or when $G$ acts diagonally on the product of its flag variety
with another variety \cite{dh}. \bl

Since $x \in X^{ss}(0)$, $G_x$ acts trivially on $(L_0)_x$.  If it
acts nontrivially on $(L_+)_x$, requiring it to act with some negative
weight $v_+ < 0$ fixes an isomorphism $G_x \cong \kst$.  It then acts
on $(L_-)_x$ with some positive weight $v_- > 0$.  To obtain the
first, stronger generalizations, assume that these two weights are
coprime: $(v_+, v_-) = 1$.  When $X$ is normal and projective and
$L_\pm$ are in adjacent chambers, this additional hypothesis can be
interpreted as follows.  The weight of the $G_x$-action defines a
homomorphism $\rho: \NS^G \to \Z$, and $L_\pm$ can be chosen within
their chambers to satisfy this hypothesis, and the conditions of the
set-up, if and only if $\rho$ is surjective.  Again, this is always
true when $G$ is a torus or when $G$ acts diagonally on the product of
its flag variety with another variety.

If $(v_+, v_-) = 1$, then $\pi^{-1}(x)$ is contained in a
$G$-orbit, so $p(\pi^{-1}(x))$ is a single point in $Z$.

\begin{s}{Lemma}
\label{1v}
If \re{1y} holds and $(v_+, v_-) = 1$, then $L_\pm$ can be chosen so
that $G$ acts freely on $Y$ at $\pi^{-1}(x)$ and $Z$ is smooth at
$p(\pi^{-1}(x))$.
\end{s}

\pf.  Since $G_x$ acts on $(L_\pm)_x$ with weights $v_\pm$, there
exist positive powers of $L_+$ and $L_-$ whose weights add to 1.
Replace $L_\pm$ by these powers.  Then $G_x$ acts freely on $(L_+
\otimes L_-^{-1})_x \sans 0$.  But this is exactly $\pi^{-1}(x)$, so
$G$ acts freely on $Y$ at $\pi^{-1}(x)$.  To show that $Z$ is smooth
at $p(\pi^{-1}(x))$, it therefore suffices to show that $\pi^{-1}(x)
\subset Y^s(G)$, that is, that the $G$-orbit of $\pi^{-1}(x)$ is
closed in $Y^{ss}(G)$.  But if $y \in Y^{ss}(G)$ is in the closure of
$G \cdot \pi^{-1}(x)$, then $y \notin i_\pm(X)$ by \re{1j}(b) and
\re{1u}(b), so $y \in \pi^{-1}(x')$ for some $x' \in X^0$ by
\re{1j}(c).  Then $x'$ is in the closure of $G \cdot x \subset X^0$,
so by \re{1y}, $x' \in G \cdot x$ and hence $y \in G \cdot
\pi^{-1}(x)$. \fp

\begin{s}{Proposition}
\label{1w}
If \re{1y} holds and $(v_+, v_-) = 1$, then \mbox{\rm (a)} $X^0$ is
smooth at $x$; \mbox{\rm (b)} on a neighbourhood of $x$ in $X^0$,
there exists a vector bundle $N$ with $\kst$-action, whose fibre at
$x$ is naturally isomorphic to $N_{G \cdot x / X}$; \mbox{\rm (c)} the
bundle $N^0$ of zero weight spaces of $N$ is exactly the image of
$TX^0$ in $N$; \mbox{\rm (d)} the bundles $N^\pm$ of positive and
negative weight spaces of $N$ are naturally isomorphic to $N_{X^0 /
X^\pm}$.
\end{s}

\pf.  By \re{2e}(a), \re{1u}(b) and \re{1v}, $Z^0$ is smooth at
$p(\pi^{-1}(x))$, and $Y^0$ is locally a principal $G$-bundle over
$Z^0$.  Hence $Y^0$ is smooth at $\pi^{-1}(x)$, so $X^0$ is smooth at
$x$.  The bundle $N_Z$ is just $TZ|_{Z^0}$, so define $N_Y = p^* N_Z$.
This is acted upon by $\kst$, so by Kempf's descent lemma \cite{dn}
descends to a bundle $N_X$ which has the desired property.  This
proves (b); the proofs of (c) and (d) are similar, using \re{2e}(c)
and (d).  \fp \bl

As in \S1, let $w_i^\pm \in \Z$ be the weights of the $\kst$-actions
on $N^\pm$.

\begin{s}{Theorem}
\label{2n}
If \re{1y} holds and $(v_+, v_-) = 1$, then over a neighbourhood of
$x$ in $X^0 \mod G(0)$, $X^\pm \mod G(\pm)$ are locally trivial
fibrations with fibre the weighted projective space $\Pj(|w_i^\pm|)$.
\end{s}

\pf.  This now follows immediately from \re{2j} and \re{1v}.  \fp \bl

If $X \mod G(\pm)$ are both nonempty, then $X^\pm \mod G(\pm)$ are the
supports of the blow-up loci of \re{1g}.  But if $X \mod G(-) =
\varnothing$, then $X^+ \mod G(+) = X \mod G(+)$ and $X^0 \mod G(0) =
X \mod G(0)$, so \re{2n} says the natural morphism $X \mod G(+) \to X
\mod G(0)$ is a locally trivial weighted projective fibration. \bl

If moreover all $w_i^\pm = \pm w$ for some $w$, then for any
linearization $L$ such that $L_x$ is acted on by $G_x$ with weight
$-1$, the bundles $N^\pm \otimes L^{\pm w}$ are acted upon trivially
by all stabilizers.  So by Kempf's descent lemma \cite{dn} they
descend to vector bundles $W^\pm$ over a neighbourhood of $x$ in $X^0
\mod G(0)$.

\begin{s}{Theorem}
\label{1x}
Suppose that \re{1y} holds, that $(v_+, v_-) = 1$, and that all
$w_i^\pm = \pm w$ for some $w$.  Then over a neighbourhood of $x$ in
$X^0 \mod G(0)$, $X^\pm \mod G(\pm)$ are naturally isomorphic to the
projective bundles $\Pj W^\pm$, their normal bundles are naturally
isomorphic to $\pi_\pm^* W^\mp(-1)$, and the blow-ups of $X \mod
G(\pm)$ at $X^\pm \mod G(\pm)$, and of $X \mod G(0)$ at $X^0 \mod
G(0)$, are all naturally isomorphic to the fibred product $X \mod G(-)
\times_{X \mod G(0)} X \mod G(+)$.
\end{s}

\pf.  First notice that, although $W^\pm$ depend on the choice of $L$,
the projectivizations $\Pj W^\pm$, and even the line bundle $\co(1,1)
\to \Pj W^+ \times_{X^0 \mod G(0)} \Pj W^-$, are independent of $L$.
Now on $Z$, taking $L \cong \co$ yields $W^\pm = N^\pm$.  But $N^\pm_Y
= p^* N^\pm_Z$, so taking $L = p^* \co$ on $Y$, with the induced
linearization, yields $W^\pm_Y = W^\pm_Z$.  On the other hand,
$N^\pm_Y = \pi^* N^\pm_X$ also, so for another choice of $L$ on $Y$,
$W^\pm_Y = W^\pm_X$.  Hence $\Pj W^\pm_X \cong \Pj W^\pm_Z$, and the
line bundles $\co(1,1) \to \Pj W^+_X \times_{X^0 \mod G(0)} \Pj W^-_X$
and $\co(1,1) \to \Pj W^+_Z \times_{Z^0} \Pj W^-_Z$ correspond under
this isomorphism; pushing down and taking duals, the bundles
$\pi_\pm^* W^\mp_X (-1) \to \Pj W^\pm_X$ and $\pi_\pm^* W^\mp_Z (-1)
\to \Pj W^\pm_Z$ also correspond.  On the other hand, by \re{1u}(a)
$X^\pm \mod G(\pm) = Z^\pm \mod T(\pm)$.  The theorem then follows
from \re{2k} together with \re{1e}.  \fp \bl

The hypothesis on $w_i^\pm$ is most easily verified as follows.

\begin{s}{Proposition}
\label{1z}
If every 0-dimensional stabilizer is trivial near $x$, then all
$w_i^\pm = \pm 1$.
\end{s}

\pf. If not all $w_i^\pm = \pm 1$, then by \re{2e}(b) there is a point
$z \in Z$ with proper nontrivial stabilizer $T_z$ such that
$p(\pi^{-1}(x))$ is in the closure of $T \cdot z$.  Then any $y \in
p^{-1}(z)$ has nontrivial 0-dimensional stabilizer $(G \times T)_y$,
and the closure of $(G \times T) \cdot y$ contains $\pi^{-1}(x)$.  But
then $G_{\pi(y)} \cong (G \times T)_y$, and the closure of $G \cdot
\pi(y)$ contains $x$.  \fp \bl

For most applications, the hypothesis \re{1y} will hold for all $x \in
X^0$.  Then the conclusions of \re{2n} and \re{1x} hold globally,
because they are natural.  Notice, however, that if $X^0$ is not
connected, then the $w_i^\pm$ need be only locally constant.

\bit{The smooth case: \'etale results}

If $(v_+, v_-) \neq 1$, however, then the proof of \re{1v} fails, and
$X^\pm \mod G(\pm)$ need not be locally trivial over $X^0 \mod G(0)$,
even if \re{1y} holds: see \re{2p} for a counterexample.  But they
will be locally trivial in the \'etale topology.  The proof uses the
Luna slice theorem.  The first step, however, is to check that $v_\pm$
are always nonzero.  As always, let $L_\pm$ and $L_0$ be as in the
set-up of \S3.

\begin{s}{Lemma}
\label{1t}
If \re{1y} holds, then $G_x$ acts nontrivially on $(L_\pm)_x$.
\end{s}

\pf.  If $G_x$ acts trivially on $L_+$ (and hence on $L_-$), then the
embedding $\kst = \pi^{-1}(x) \subset Y$ descends to an embedding
$\kst \subset Z$.  This is completed by two points, which must come
from two equivalence classes of $G$-orbits in $Y^{ss}(G)$.  These
semistable orbits cannot be $G \cdot (i_\pm(x))$, since $i_\pm(x)
\notin Y^{ss}(G)$ by \re{1j}(a) and the definition of $X^0$.  Hence
our two classes of semistable $G$-orbits must be contained in
$\pi^{-1}(\overline{G \cdot x} \sans G \cdot x)$.  By \re{1j}(b) and
(c) their images in $\pi$ are in $X^{ss}(0)$.  But they are also in
the closure of $G \cdot x$, which contradicts \re{1y}.  \fp

Again, requiring $G_x$ to act on $(L_+)_x$ with negative weight $v_+ <
0$ fixes an isomorphism $G_x \cong \kst$.  It then acts on $(L_-)_x$
with positive weight $v_- > 0$.  We no longer require $(v_+, v_-) =
1$, but assume instead the following.

\begin{r}{Hypothesis}
\label{2s}
Suppose that either $\chr \k = 0$ or $(v_+, v_-)$ is coprime to
$\chr \k$.
\end{r}

Now choose $y \in \pi^{-1}(x)$, and let $S = (G \times T)_y$.

\begin{s}{Lemma}
\label{1o}
If \re{1y} holds, then there is a fixed isomorphism $S \cong
\kst$, and $G_y = S \cap G$ is a proper subgroup such that $S / G_y
\cong T$ naturally.  Moreover, if \re{2s} holds, then $L_\pm$ may be
chosen so that $G_y$ is reduced.
\end{s}

\pf.  Since $\pi$ is the quotient morphism for the $T$-action, $S
\subset G_x \times T = \kst \times T$; indeed, it is the subgroup
acting trivially on $\pi^{-1}(x)$.  Since by \re{1t} $\kst$ acts
nontrivially, and $T$ obviously acts with weight 1, this has a fixed
isomorphism to $\kst$, and its intersection with $\kst \times 1$ is a
proper subgroup having the desired property.

If $\chr \k = 0$, then $G_y$ is certainly reduced.  Otherwise, let $v
= (v_+, v_-)$, and replace $L_+$ and $L_-$ with positive powers so
that $G_x \cong \kst$ acts on $\pi^{-1}(x) = (L_+ \otimes L_-^{-1})_x$
with weight $v$.  Then $G_y = \Spec \k[z]/ \langle z^v - 1 \rangle$,
which is reduced if $v$ is coprime to $\chr \k$.  \fp \bl

Notice that for any $y \in \pi^{-1}(x)$, $(G \times T) \cdot y = G
\cdot y$.  Together with lemma \re{1o}, this implies that any
$S$-invariant complement to $T_y (G \times T) \cdot y$ in $T_y Y$ is
also a $G_y$-invariant complement to $T_y (G \cdot y)$ in $T_y Y$.  It
follows from the definition of the Luna slice \cite{luna,mf} that a
slice for the $(G \times T)$-action at $y$ is also a slice for the
$G$-action at $y$.  Luna's theorem then implies that there exists a
smooth affine $U \subset Y$ containing $y$ and preserved by $S$, and a
natural diagram
$$\begin{array}{ccccc}
     G \times U & \lrow  & G
\times_{G_y} U  & \lrow & Y^{ss}(G) \vspace{.7ex} \\
     \down{} &  & \down{} &  & \down{} \vspace{.7ex} \\
     U & \lrow & U / G_y & \lrow & Z,
\end{array} $$
such that the two horizontal arrows on the right are strongly \'etale
with respect to the actions of
$$\begin{array}{ccccc}
     G \times S &\lrow & G \times T  & \lrow  &  G \times T \vspace{.7ex} \\
     \down{} &  & \down{} &  & \down{} \vspace{.7ex} \\
     S &\lrow & T  & \lrow  & T.
\end{array} $$

These actions are obvious in every case except perhaps on $G \times
U$; there the $G \times S$-action is given by $(g, s)\cdot(h,u)
= (gh\hat{s}^{-1}, s u)$, where $\hat{s}$ is the
image of $s$ in the projection $S \to G$.  Each of these actions
has a 1-parameter family of fractional linearizations, pulled back from the
right-hand column.  For any object $V$ in the diagram, define $V^0$
and $V^\pm$ with respect to these linearizations.

\begin{s}{Lemma}
\label{1k}
For every arrow $f: V \to W$ in the diagram, $f^{-1}(W^\pm) = V^\pm$
and $f^{-1}(W^0) = V^0$.
\end{s}

\pf.  This is straightforward for the vertical arrows, and for the
morphism $U \to U / G_y$, because they are all quotients by subgroups
of the groups which act.  The result for $G \times U \lrow  G
\times_{G_y} U$ follows from the commutativity of the diagram.  As for
the strongly \'etale morphisms, these are treated as in the proof of
\re{2e}(b).  \fp \bl

With this construction, \re{1w} can now be strengthened.

\begin{s}{Proposition}
If \re{1y} and \re{2s} hold, then the conclusions of \re{1w} follow
even if $(v_+, v_-) \neq 1$.
\end{s}

\pf.  By \re{1o} and \re{1k}, in a neighbourhood of $\pi^{-1}(x)$ in
$Y$, a point is in $Y^0$ if and only if it has stabilizer conjugate to
$S$.  Hence in a neighbourhood of $x$ in $X$, a point is in $X^0$ if
and only if it has stabilizer conjugate to $G_x$.  So if $W \subset
X^0$ is the fixed-point set for the $G_x$-action on a neighbourhood of
$x \in X$, then a neighbourhood of $x$ in $X^0$ is precisely the
affine quotient of $W \times G$ by the diagonal action of the
normalizer of $G_x$ in $G$.  In particular, this is smooth as claimed
in (a), because $W$ is smooth by \re{2e}(a) and the normalizer acts
freely.  Also, its tangent bundle has a natural subbundle consisting
of the tangent spaces to the $G$-orbits.  Let $N_X$ be the quotient of
$TX|_{X^0}$ by this subbundle.  Then $N_X$ certainly satisfies (b),
and $N_Y = \pi^* N_X$, so by it suffices to prove (c) and (d) for $Y$
with its $G \times T$-action.  But (c) and (d) hold for $U$ with its
$S$-action by \re{2e}(c) and (d), hence for $G \times U$ with its $G
\times S$-action since $N_U$ pulls back to $N_{G \times U}$.  But the
morphism $G \times_{G_y} U \to Y$ is \'etale, and so is the morphism
$G \times U \to G \times_{G_y} U$, since $G_y$ is reduced.  The result
for $Y$ then follows from \re{1k}, since \'etale morphisms are
isomorphisms on tangent spaces.  \fp \bl

Let $w_i^\pm \in \Z$ be the weights of the $\kst$-action on $N^\pm$.

\begin{s}{Theorem}
\label{1l}
If \re{1y} and \re{2s} hold, then over a neighbourhood of $x$ in $X^0
\mod G(0)$, $X^\pm \mod G(\pm)$ are fibrations, locally trivial in the
\'etale topology, with fibre the weighted projective space
$\Pj(|w_i^\pm|)$.
\end{s}

As before, if $X \mod G(\pm)$ are both nonempty, then $X^\pm \mod
G(\pm)$ are the supports of the blow-up loci of \re{1g}; but if $X
\mod G(-) = \varnothing$, then \re{1l} says the natural morphism $X
\mod G(+) \to X \mod G(0)$ is a weighted projective fibration, locally
trivial in the \'etale topology.

The proof requires the following lemma.

\begin{s}{Lemma}
\label{1m}
If $\phi: V \to W$ is a strongly \'etale morphism of affine varieties
with $\kst$-action, then $\phi \mod \pm : V \mod \pm \to W \mod \pm$
are \'etale, and $V \mod \pm = {W \mod \pm} \times_{W \mod 0} {V \mod
0}$.
\end{s}

\pf.  Say $V = \Spec R$, $W = \Spec S$.  The $\kst$-actions induce
$\Z$-gradings on $R$ and $S$, and $V = W \times_{W \mod 0} V \mod 0$
implies $R = S \otimes_{S_0} R_0$.  Hence $\bigoplus_{i \in \N} R_{\pm
i} = \bigoplus_{i \in \N} S_{\pm i} \otimes_{S_0} R_0$, which implies
the second statement.  Then $\phi \mod \pm$ are certainly \'etale,
since being \'etale is preserved by base change.  \fp

\pf\ of \re{1l}.  The Luna slice $U$ associated to any $y \in
\pi^{-1}(x)$ is smooth, and for $S = (G \times T)_y$, $U \mod S(t) =
(U / G_y) \mod T(t)$ since $T = S / G_y$ by \re{1o}.  But as stated
when $U$ was constructed, $U / G_y$ is strongly \'etale over $Z$, so
by \re{1m} and \re{1e} $U \mod S(t)$ is \'etale over $Z \mod T(t) = X
\mod G(t)$, and $U \mod S(\pm) = X \mod G(\pm) \times_{X \mod G(0)} U
\mod S(0)$.  In particular, $U^0 \mod S(\pm) = X^0 \mod G(\pm)
\times_{X^0 \mod G(0)} U^0 \mod S(0)$, which is exactly the pullback
of $X^0 \mod G(\pm)$ by the \'etale morphism $U^0 \mod S(0) \to X^0
\mod G(0)$.  The theorem therefore follows from the analogous result
\re{2j} for quotients of smooth affines by $\kst$.  \fp

\begin{s}{Counterexample}
\label{1r}
To show that the hypothesis $G_x \cong \kst$ is necessary in \re{1l}.
\end{s}

Let $G$ be any semisimple reductive group, and let $V_+$ and $V_-$ be
representations of $G$.  Let $X = \Pj (V_+ \oplus V_- \oplus \k)$, and
let $G \times \kst$ act on $X$, $G$ in the obvious way, and $\kst$
with weights $1,-1,0$.  Then $\NS_\Q^{G \times \kst} \cong \Q$, with
two chambers separated by a wall at 0.  Moreover $X^\pm = \Pj(V_\pm
\oplus \k)$, so $X^0 = \{ (0,0,1) \}$.  But $G_{(0,0,1)} = G \times
\kst$, so the hypothesis is violated.  Now $X^\pm \mod \kst (\pm) =
\Pj (V_\pm)$, so $X^\pm \mod (G \times \kst)(\pm) = \Pj (V_\pm) \mod
G$.  This certainly need not be a projective space, as the theorem
would predict; see for example the discussion of the case $G = {\rm
PSL(2)}$ in \S6.  \fp \bl

Since $G$ acts on $N^\pm$, there are quotients $N^\pm \mod G(\pm)$,
which are fibrations with fibre $\Pj(|w_i^\pm|)$ over a neighbourhood
of $x$ in $X^0 \mod G(0)$, locally trivial in the \'etale topology.
Notice that by \re{1k}, since $N^\pm_V = N_{V^0 / V^\pm}$ for $V = X$,
$Y$, $G \times_{G_y} U$, and $U$,
\beqas
N^\pm_U \mod S(\pm) & = & N^\pm_{G \times U} \mod (G \times S)(\pm) \\
 & = & N^\pm_Y \mod (G \times T)(\pm) \times_{Y^0 \mod (G \times
T)(0)} (G \times U)^0 \mod (G \times S)(0) \\
 & = & N^\pm_X \mod G(\pm) \times_{X^0 \mod G(0)} U^0 \mod
S(0),
\eeqas
which is the pullback of $N^\pm_X \mod G(\pm)$ by the \'etale morphism
$U^0 \mod S(0) \to X^0 \mod G(0)$.  The following result then ought to
be true, but proving it conclusively is rather cumbersome, so we content
ourselves with a sketch.

\begin{s}{Theorem}
\label{1n}
Suppose that \re{1y} and \re{2s} hold, and that all $w_i^\pm = \pm w$
for some $w$.  Then $X^\pm \mod G(\pm)$ are naturally isomorphic to
$N^\pm \mod G(\pm)$, and the blow-ups of $X \mod G(\pm)$ at $X^\pm
\mod G(\pm)$, and of $X \mod G(0)$ at $X^0 \mod G(0)$, are all
naturally isomorphic to the fibred product $X \mod G(-) \times_{X \mod
G(0)} X \mod G(+)$.
\end{s}

{\em Sketch of proof}.  All the blow-ups and the fibred product are
empty if either $X \mod G(+)$ or $X \mod G(-)$ is empty, so suppose
they are not.  Now $X^\pm \mod G(\pm) \to X^0 \mod G(0)$ are covered
in the \'etale topology by $U^\pm \mod S(\pm) \to U^0 \mod S(0)$ by
\re{1m}, and $N^\pm_X \mod G(\pm) \to X^0 \mod G(0)$ are covered in
the \'etale topology by $N^\pm_U \mod S(\pm) \to U^0 \mod S(0)$ by the
remarks above.  But the analogous result for $U$ holds by \re{2k}.
The theorem would therefore follow if we could display a morphism
$N^\pm_X \mod G(\pm) \to X^\pm \mod G(\pm)$ compatible with the
\'etale morphisms and the isomorphisms of \re{2k}.  Unfortunately,
this is somewhat awkward to construct.  One way to do it is to imitate
the argument of \re{1x}, using a bundle of tangent cones with
$\kst$-action over $Z^0$, which is typically in the singular locus of
$Z$.  This requires generalizing the Bialynicki-Birula decomposition
theorem to the mildly singular space $Z$, which can still be
accomplished using the Luna slice theorem.  \fp \bl

The hypothesis on $w_i^\pm$ is again most easily verified as follows.

\begin{s}{Proposition}
If every 0-dimensional stabilizer is trivial near $x$, then all
$w_i^\pm = \pm 1$.
\end{s}

\pf. If not all $w_i^\pm = \pm 1$, then for each Luna slice $U$ there
is a point $u \in U$ with nontrivial proper stabilizer $S_u$.  Then
any $(g,u) \in G \times U$ satisfies $(G \times S)_{(g,u)} \cong S_u$.
Since the morphism $G \times U \to Y$ is \'etale, this implies that
there exists $y \in Y$ with a nontrivial 0-dimensional stabilizer $(G
\times T)_y$.  But then $G_{\pi(y)} = (G \times T)_y$. \fp \bl

Again, for most applications, the hypothesis \re{1y}, and hence the
conclusions of \re{1l} and \re{1n}, will hold globally.

\bit{The first example}

In this section we turn to a simple application of our main results,
the much-studied diagonal action of $\PSL{2}$ on the $n$-fold product
$(\Pj^1)^n$.  This has $n$ independent line bundles, so it is tempting
to study the quotient with respect to an arbitrary $\co(t_1, \dots,
t_n)$.  We will take a different approach, however: to add an $n+1$th
copy of $\Pj^1$, and consider only fractional linearizations on
$(\Pj^1)^{n+1}$ of the form $\co(t,1,1,\dots,1)$.  This has the
advantage that it does not break the symmetry among the $n$ factors.
In other words, the symmetric group $S_n$ acts compatibly on
everything, so in addition to $(\Pj^1)^n$, we learn about quotients by
$\PSL{2}$ of the symmetric product $(\Pj^1)^n / S_n = \Pj^n$.

So for any $n>2$, let $(\Pj^1)^n$ be acted on diagonally by $G =
\PSL{2}$, fractionally linearized on $\co(1,1,\dots,1)$.  We wish to
study the quotient $(\Pj^1)^n \mod G$.  The stability condition for
this action is worked out in \cite[4.16; GIT Ch.\ 3]{n}, using the
numerical criterion.  This is readily generalized to an arbitrary
linearization on $X = \Pj^1 \times (\Pj^1)^n = (\Pj^1)^{n+1}$; indeed
for the fractional linearization $\co(t_0, \dots, t_n)$, it turns out
that $(x_j) \in X^{ss}(t_j)$ if and only if, for all $x \in \Pj^1$,
$$\sum_{j=0}^n t_j \, \delta(x,x_j) \leq \sum_{j=0}^n t_j / 2.$$
Moreover, $(x_j) \in X^s(t_j)$ if and only if the inequality is always
strict.  We will study the case where $t_0$ is arbitrary, but $t_j =
1$ for $j>0$.

For $t_0 < 1$, it is easy to see that $\Pj^1 \times ((\Pj^1)^n)^s
\subset X^{ss} \subset \Pj^1 \times ((\Pj^1)^n)^{ss}$.  So the
projection $X \to (\Pj^1)^n$ induces a morphism $X \mod G(t_0) \to
(\Pj^1)^n \mod G$ whose fibre over each stable point is $\Pj^1$.
Indeed, each diagonal in $X = \Pj^1 \times (\Pj^1)^n$ is fixed by $G$,
so descends to $X \mod G(t_0)$.  Hence, over the stable set in
$(\Pj^1)^n \mod G$, $X \mod G(t_0)$ is exactly the total space of the
universal family.

Now because $G = \PSL{2}$, not $\SL{2}$, the bundle
$\co(1,0,0,\dots,0)$ has no bona fide linearization, only a fractional
one.  However, $\co(1,1,0,\dots,0)$ does admit a bona fide
linearization, as does $\co(1,1,\dots,1)$ if $n$ is odd.  So these
bundles descend to $X \mod G(t_0)$ for $t_0 < 1$, yielding line
bundles whose restriction to each $\Pj^1$ fibre is $\co(1)$.  This
implies that, over the stable set in $(\Pj^1)^n \mod G$, $X \mod
G(t_0)$ is a locally trivial fibration.  In particular, if $n$ is odd,
it is a fibration everywhere.  However, if $n$ is even, there is no
$S_n$-invariant line bundle having the desired property.  Hence the
quotient $(\Pj^1 \times \Pj^n) \mod G(t_0) = (X \mod G(t_0)) / S_n$,
though it has generic fibre $\Pj^1$ over $\Pj^n \mod G$, and is
generically trivial in the \'etale topology, is not even locally
trivial anywhere.  It is (generically) what is sometimes called a {\em
conic bundle}.

To apply our results to this situation, note first that for numerical
reasons the stability condition only changes when equality can occur
in the inequality above, that is, when $t_0 = n - 2m$ for some integer
$m \leq n/2$.  These will be our walls.  A point $(x_j) \in X$ is in
$X^0$ for one of these walls if it is semistable for $t = t_0$, but
unstable otherwise.  This means there exist points $x, x' \in \Pj^1$
such that
$$t \, \delta(x, x_0) + \sum_{j=1}^n \delta(x,x_j) \geq \frac{t+n}{2} $$
for $t \leq t_0$, and
$$t \, \delta(x', x_0) + \sum_{j=1}^n \delta(x',x_j) \geq \frac{t+n}{2} $$
for $t \geq t_0$, with equality in both if and only if $t = t_0$.
This requires that all $x_j$ be either $x$ or $x'$; indeed, $x_0$ and
exactly $m$ other $x_j$ must be $x$, and the $n-m$ remaining $x_j$
must be $x'$.  In particular, this implies $m \geq 0$, so there are
only finitely many walls, as expected. On the other hand, any $(x_j)
\in X$ of this form will belong to $X^{ss}(t_0)$, provided that $x
\neq x'$.  So $X^0(t_0)$ consists of ${n \choose m}$ copies of $(\Pj^1
\times \Pj^1) \sans \Delta$.

Hence every point $(x_j) \in X^0(t_0)$ is stabilized by the subgroup
of $G$ fixing $x$ and $x'$, which is isomorphic to $\kst$.  So the
hypothesis \re{1y} is satisfied.  Moreover, the bundle
$\co(2,0,0,\dots,0)$ is acted on by this $\kst$ with weight 1, so the
strong results of \S4 will apply.  Finally, we claim that, even though
$X^0$ is disconnected, the weights $w_i^\pm$ are globally constant,
and are all $\pm 1$.  Indeed, it is easy to see that the $w_i^\pm$ are
independent of the component, because the action of the symmetric
group $S_n$ on $X$ commutes with the $G$-action, and acts transitively
on the components of $X^0$.  To evaluate $w_i^\pm$, note that each
component is a single orbit, and that setting $x=0$, $x' = \infty$
determines an unique point in this orbit with stabilizer $\{ \diag
(\lambda^{-1}, \lambda) \st \lambda \in \kst \} / \pm 1$.  This acts
on $T\Pj^1$ with weight $-1$ at 0, $1$ at $\infty$; so it acts on
$T_{(x_j)}X$ with $m+1$ weights equal to $-1$ and $n-m$ weights equal
to $1$. But it acts on the $G$-orbit $G / \kst$ with one weight equal
to 1 and one equal to $-1$, so $N$ is acted on with $m$ weights $-1$
and $n-m-1$ weights $1$.  So the very strongest result \re{1x}
applies.  Hence $X^\pm \mod G(\pm)$ are bundles with fibre
$\Pj^{n-m-2}$ and $\Pj^{m-1}$, respectively, over $X^0 \mod G(0)$.
Since this is just ${n \choose m}$ points, $X^\pm \mod G(\pm)$ are
disjoint unions of ${n \choose m}$ projective spaces.  Moreover, the
blow-ups of $X \mod G(\pm)$ at $X^\pm \mod G(\pm)$ are both isomorphic
to $X \mod G(-) \times_{X \mod G(0)} X \mod G(+)$.

This does not seem to say much about $(\Pj^1)^n \mod G$ itself, only
about $X \mod G(t)$, which for $t$ small is (at least generically) a
$\Pj^1$-bundle over it.  But this is enough to compute quite a lot
(cf.\ \cite{t1}).  We content ourselves with just one calculation, of
the Betti numbers of $(\Pj^1)^n \mod G$ and $\Pj^n \mod G$ for $n$
odd, in the case where the ground field is the complex numbers $\C$.
These formulas are originally due to Kirwan \cite{k}.

\begin{s}{Proposition}
For $n$ odd,
$$P_t ((\Pj^1)^n \mod G) = \sum_{m=0}^{(n-1)/2} {n \choose m} \frac{t^{2m} -
t^{2(n-m-1)}}{1-t^4} $$
and
$$P_t (\Pj^n \mod G) = \sum_{m=0}^{(n-1)/2} \frac{t^{2m} -
t^{2(n-m-1)}}{1-t^4}. $$
\end{s}

\pf.  Let $t_0 = n - 2m$, and $t_\pm = t_0 \mp 1$.  Then the blow-ups
of $X \mod G(\pm)$ at $X^\pm \mod G(\pm)$ are equal.  So by the
standard formula for Poincar\'e polynomials of blow-ups,
$$ P_t(X \mod G(-)) - P_t(X^- \mod G(-)) + P_t(E)
=  P_t(X \mod G(+)) - P_t(X^+ \mod G(+)) + P_t(E), $$
where $E$ is the exceptional divisor.  Cancelling and rearranging
yields
$$ P_t(X \mod G(+)) - P_t(X \mod G(-))
=  P_t(X^+ \mod G(+)) - P_t(X^- \mod G(-)). $$
But $X^\pm \mod G(\pm)$ are ${n \choose m}$ copies of $\Pj^{n-m-2}$
and $\Pj^{m-1}$ respectively, so
\beqas
P_t(X \mod G(+)) - P_t(X \mod G(-))
& = & {n \choose m}
\left( \frac{1-t^{2(n-m-1)}}{1-t^2}-\frac{1-t^{2m}}{1-t^2} \right) \\
& = & {n \choose m} \frac{t^{2m}-t^{2(n-m-1)}}{1-t^2}.
\eeqas
Summed over $m$, the left-hand side telescopes, so for $t<1$
$$P_t(X \mod G(t)) = \sum_{m=0}^{(n-1)/2} {n \choose m}
\frac{t^{2m}-t^{2(n-m-1)}}{1-t^2}. $$
But for $n$ odd, this is a $\Pj^1$-bundle over $(\Pj^1)^n \mod G$, and the
Poincar\'e polynomial of any projective bundle splits, so
$$P_t((\Pj^1)^n \mod G) = \sum_{m=0}^{(n-1)/2} {n \choose m}
\frac{t^{2m}-t^{2(n-m-1)}}{1-t^4}, $$
as desired.

As for $\Pj^n \mod G$, it is the quotient of $(\Pj^1)^n \mod G$ by the
action of the symmetric group $S_n$.  A result from Grothendieck's
T\^ohoku paper then implies \cite{gr,mac} that $H^*(\Pj^n \mod G, \C)$
is the $S_n$-invariant part of $H^*((\Pj^1)^n \mod G, \C)$.  But since
$S_n$ acts on $X \mod G(t)$ for all $t$, the calculation above
actually decomposes $H^*((\Pj^1)^n \mod G, \C)$ as a representation of
$S_n$: the term with coefficient ${n \choose m}$ gives the
multiplicity of the permutation representation induced by the natural
action of $S_n$ on subsets of $\{ 1, \dots, n\}$ of size $m$.  The
trivial summand of this representation is exactly one-dimensional, so
the cohomology of $\Pj^n \mod G$ is as stated.  \fp \bl

We round off this section by using some of the ideas discussed above
to give the counterexample promised in \S5.

\begin{s}{Counterexample}
\label{2p}
To show that the hypothesis $(v_+, v_-) = 1$ is necessary in \re{2n}.
\end{s}

Let $V$ be the standard representation of $\GL{2}$, and let $W = S^n V
\otimes (\Lambda^2 V)^{-n/2}$, where $S^n V$ is the $n$th symmetric
power for some even $n>2$.  Let $X = \Pj (V \oplus V^* \oplus W)$, and
let $\GL{2}$ act on $X$.  Then $\NS_\Q^{\GL{2}} \cong \Q$, with two
chambers separated by a wall at 0.  The central $\kst \subset \GL{2}$
acts on $V$, $V^*$, and $W$ with weight $1$, $-1$, and $0$
respectively, so $X^+$ is open in $\Pj(V \oplus W)$, $X^-$ is open in
$\Pj(V^* \oplus W)$, and $X^0$ is open in $\Pj W$.  By construction,
the scalars $\kst \subset \GL{2}$ act trivially on $\Pj W$ with the
linearization $L_0$, so the action reduces to the action of $\PSL{2}$
on $\Pj^n$ considered above.  A generic $x \in \Pj^n$ is stable and is
acted on freely by $\PSL{2}$, so $\GL{2}_x = \kst$.  Moreover, it is
stable, so $\GL{2} \cdot x = \PSL{2} \cdot x$ is closed in
$(\Pj^n)^{ss}$ and hence in $X^{ss}(0)$.  Therefore \re{1y} holds for
the $\GL{2}$-action at $x$.  On the other hand, the hypothesis $(v_+,
v_-) = 1$ cannot be satisfied: the tautological linearization $L_0$ on
$\co(1)$ is acted on with weight 0, and the linearization $L_+$
obtained by tensoring $L_0$ with the character $\det: \GL{2} \to \kst$
is acted on with weight $2$, but together these generate $\Pic^\GL{2}
X$.

Now $X^+ \mod \kst (+) = \Pj V \times \Pj W$.  As a variety with
$\PSL{2}$-action, this is exactly $\Pj^1 \times \Pj^n$ as considered
above.  The $+$ linearization on $\Pj V \times \Pj W$ corresponds to
the linearization given by $t<1$ on $\Pj^1 \times \Pj^n$, so $X^\pm
\mod \GL{2}(+) = (\Pj V \times \Pj W) \mod \PSL{2}(+) = (\Pj^1 \times
\Pj^n) \mod \PSL{2}(t)$.  As mentioned above, this is a conic bundle,
so it is not even locally trivial over $X^0 \mod \GL{2}(0) = \Pj W
\mod \PSL{2}(0) = \Pj^n \mod \PSL{2}$ at $x$.  \fp

\bit{Parabolic bundles}

In the last two sections we apply our main results to moduli problems
of vector bundles with additional structure over a curve.  Throughout
these sections, $C$ will denote a smooth projective curve over $k$, of
genus $g>0$.

Fix a point $p \in C$.  In this section we will study parabolic
bundles of rank $r$ and degree $d$ over $C$, with parabolic structure
at $p$.  We refer to \cite{ms,s} for basic definitions and results on
parabolic bundles.  However, we insist for simplicity on full flags at
$p$, so the weights $\ell_j \in [0,1)$ are strictly increasing.  The
space of all possible weights is therefore parametrized by $$W = \{
(\ell_j) \in \Q^r \st 0 \leq \ell_1 < \ell_2 < \cdots < \ell_r < 1
\}.$$

There are several constructions of the moduli space $M(\ell_j)$ of
parabolic bundles semistable with respect to $(\ell_j)$.  The one best
suited for our purposes is due to Bhosle \cite{bho}, following
Gieseker \cite{g}; so we first review his construction, then hers.

Suppose without loss of generality that $d > > 0$, and let $\chi = d +
r(1-g)$.  Let $\Quot$ be the Grothendieck Quot scheme \cite{quot}
parametrizing quotients $\phi: \co_C^\chi \to E$, where $E$ has
Hilbert polynomial $\chi +ri$ in $i$, and let $\co^\chi \to \be$ be
the universal quotient over $\Quot \times C$.  Let $R \subset \Quot$
be the smooth open subvariety consisting of locally free sheaves $E$
such that $H^0(\co^\chi) \to H^0(E)$ is an isomorphism, and let
$R^{ss}$ be the subset corresponding to semistable bundles.  For $d$
large, every semistable bundle of rank $r$ and degree $d$ is
represented by a point in $R$.  Let $Z$ be the bundle over $\Pic^d C$,
constructed as a direct image, with fibre $\Pj \Hom (H^0(\La^r
\co^\chi), H^0(M))$ at $M$.  The group $G = \SL{\chi}$ acts on $R$ and
$Z$, and there is a natural $G$-morphism $T: R \to Z$, and a
linearization $L$ on $Z$, such that $T^{-1}Z^{ss}(L) = R^{ss}$.
Moreover, the restriction $T: R^{ss} \to Z^{ss}(L)$ is finite.  The
existence of a good quotient of $R^{ss}$ by $G$ then follows from a
lemma \cite[Lemma 4.6]{g} which states that if a set has a good
$G$-quotient, then so does its preimage by a finite $G$-morphism.
This quotient is the moduli space of semistable bundles on $C$.

To construct the moduli space of semistable parabolic bundles in an
analogous way, let $\tilde{R}$ be the bundle $\Fl \be|_{R \times \{ p
\} }$ of full flags in $\be_p$. This parametrizes a family of
quasi-parabolic bundles; for $d$ large, any bundle which is semistable
for some weights $(\ell_j)$ is represented by a point in $\tilde{R}$.
Let $\tilde{R}^{ss}(\ell_j)$ be the subset corresponding to parabolic
bundles semistable with respect to $(\ell_j)$.  Also let $G_r$ be the
product of Grassmannians $\prod_{j=1}^r \Gr(\chi - j,\chi)$.  Then $G$
acts on $\tilde{R}$, $Z$, and $G_r$, and there is a $G$-morphism
$\tilde{T}: \tilde{R} \to Z \times G_r$, and a family $L(\ell_j)$ of
fractional linearizations on $Z \times G_r$ depending affinely on
$(\ell_j)$, such that $\tilde{T}^{-1}(Z \times G_r)^{ss}(\ell_j) =
\tilde{R}^{ss}(\ell_j)$.  Moreover, the restriction $\tilde{T}:
\tilde{R}^{ss}(\ell_j) \to (Z \times G_r)^{ss}(\ell_j)$ is finite.
The existence of a good quotient $M(\ell_j)$ again follows from the
lemma.

In fact, we can say more.

\begin{s}{Proposition}
\label{2u}
$T$ is an embedding.
\end{s}

\pf. Since $T$ is injective \cite[4.3]{g}, it suffices to show its
derivative is everywhere injective.  At a quotient $\phi: \co^\chi \to
E$, the tangent space to $\Quot$ is given by the hypercohomology group
$\Hyp^1(\End E \stackrel{\phi}{\to} E \otimes \co^\chi)$ (cf.\
\cite{bdw,t1}).  Since $H^1(E \otimes \co^\chi) = \k^\chi \otimes
H^1(E) = 0$, this surjects onto $H^1(\End E)$, and hence onto
$H^1(\co)$, which is the tangent space to $\Pic^d C$.  So it suffices
to show the kernel of this surjection injects into the tangent space
to $\Pj \Hom(H^0(\La^r \co^\chi), H^0(\La^r E))$.  The kernel is
isomorphic to the quotient of $\Hyp^1(\End_0 E \stackrel{\phi}{\to} E
\otimes \co^\chi)$ by the 1-dimensional subspace generated by $\phi$
(cf.\ \cite[2.1]{t1}), where $\End_0$ denotes trace-free
endomorphisms.  So the kernel injects as desired if and only if the
natural map $\Hyp^1(\End E \stackrel{\phi}{\to} E \otimes \co^\chi)
\to \Hom(H^0(\La^r \co^\chi), H^0(\La^r E))$ is injective as well.

What is this natural map?  It is obtained from the derivative of $T$;
since $T$ is essentially $\phi \mapsto \La^r \phi$, the derivative of
$T$ at $\phi$ is essentially $\psi \mapsto (\La^{r-1} \phi) \wedge
\psi$.  More precisely, an element of the hypercohomology group above
is determined by \v Cech cochains $g \in C^1(\End_0 E)$ and $\psi
\in C^0(E \otimes \co^\chi)$ such that $g \phi = d\psi$. Since $\phi$
is surjective, the hypercohomology class of the pair is uniquely
determined by $\psi$; on the other hand, a cochain $\psi$ determines
the trivial hypercohomology class if and only if $\psi = f \phi$ for
some $f \in C^0(\End_0 E)$.  The natural map to $\Hom(H^0(\La^r
\co^\chi), H^0(\La^r E))$ is then indeed given by $\psi \mapsto
(\La^{r-1} \phi) \wedge \psi$; its injectivity follows from the lemma
below.  \fp

\begin{s}{Lemma}
If $\phi: \k^\chi \to \k^r$ is a linear surjection of vector spaces,
and $\psi: \k^\chi \to \k^r$ is a linear map, then $(\La^{r-1} \phi)
\wedge \psi = 0$ if and only if $\psi = f \phi$ for some $f \in \End_0
\k^r$.
\end{s}

\pf.  Suppose first that $\chi = r$.  Then $\phi$ is invertible, and
$(\La^{r-1} \phi) \wedge \psi$ is a homomorphism of 1-dimensional
vector spaces.  Indeed, if $e_1, \dots, e_r$ is the standard basis for
$\k^r$, then $\La^r \k^\chi$ is spanned by $\phi^{-1}e_1 \wedge \cdots
\wedge \phi^{-1}e_r$.  But
\beqas
((\La^{r-1} \phi) \wedge \psi)(\phi^{-1}e_1 \wedge \cdots \wedge
\phi^{-1}e_r) & \!\!\!\!\!\!\!= & \!\!\!\!\!\!\!1/r \sum_i e_1 \wedge \cdots
\wedge e_{i-1} \wedge
\psi\phi^{-1}(e_i) \wedge e_{i+1} \wedge \cdots \wedge e_r \\
& \!\!\!\!\!\!\! = & \!\!\!\!\!\!\! 1/r \, (\tr \psi\phi^{-1})(e_1 \wedge
\cdots \wedge e_r),
\eeqas
so the lemma is true when $\chi = r$.  In the general case,
$(\La^{r-1} \phi) \wedge \psi = 0$ implies $\ker \psi \supset \ker
\phi$, for if not, let $\phi^{-1}$ be a right inverse for $\phi$, and
let $u \in \ker \phi \sans \ker \psi$.  For some $i$, the coefficient
of $e_i$ in $\psi(u)$ is nonzero.  Then
\beqas
((\La^{r-1} \phi) \wedge \psi) (\phi^{-1}e_1 \wedge \cdots \wedge
\phi^{-1}e_{i-1} \wedge u \wedge \phi^{-1}e_{i+1} \wedge \cdots \wedge
\phi^{-1}e_r) & & \\
& \hspace{-60ex} = & \hspace{-31ex} e_i \wedge \cdots \wedge e_{i-1} \wedge
\psi(u) \wedge
e_{i+1} \wedge \cdots \wedge e_r \neq 0.
\eeqas
So $\psi$ descends to $\k^\chi / \ker \phi$; this has dimension $r$,
so the case above applies.  \fp

\begin{s}{Corollary}
\label{2t}
$\tilde{T}$ is an embedding.
\end{s}

\pf. $\tilde{R}$ is a bundle of flag varieties over $R$, and each
fibre clearly embeds in $G_r$. \fp

Let $X$ be the Zariski closure in $Z \times G_r$ of $\tilde{T}
(\tilde{R})$.

\begin{s}{Corollary}
The moduli space $M(\ell_j)$ of semistable parabolic bundles is $X
\mod G(\ell_j)$.
\end{s}

\pf.  Since $\tilde{T}^{-1}(Z \times G_r)^{ss}(\ell_j) =
\tilde{R}^{ss}(\ell_j)$, this is automatic provided there are no
semistable points in $X \sans \tilde{T}(\tilde{R})$.  Since
$M(\ell_j)$ is already projective, any such points would be in the
orbit closures of semistable points in $\tilde{T}(\tilde{R})$.  Hence
there would be $x \in \tilde{R}^{ss}(\ell_j)$, and a 1-parameter
subgroup $\la(t) \subset G$, such that $\lim \la(t) \cdot x \notin
\tilde{T}(\tilde{R})$, but $\mu^{\ell_j} (\tilde{T}(x), \la) = 0$,
where $\mu^{\ell_j}$ is the valuation used in the numerical criterion
\cite[4.8; GIT Defn.\ 2.2]{n}.  But all the destabilizing subgroups of
points in $\tilde{T}(\tilde{R})$ correspond to destabilizing
subbundles, and their limits are points corresponding to the
associated graded subbundles; in particular, they are in
$\tilde{T}(\tilde{R})$.  \fp

We are therefore in a position to apply our main results.  Let us
first look for walls and chambers.  Notice that the stability
condition only changes at values where there can exist subbundles
whose parabolic slope equals that of $E$.  If such a subbundle has
rank $r^+$, degree $d^+$, and weights $\ell_{j^+_i}$ for some $ \{
j^+_i \} \subset \{ 1, \dots, r \}$, then the slope condition is $$
\frac{d^+ + \sum_{i = 1}^{r^+} \ell_{j^+_i}}{r^+} = \frac{d + \sum_{j
= 1}^r \ell_j}{r}. $$ This determines a codimension 1 affine subset of
$W$, which is one of our walls.  The complementary numbers $r^- = r -
r^+$, $d^- = d - d^+$, and $\{ \ell_{j^-_i} \} = \{ \ell_j \} \sans \{
\ell_{j^+_i} \}$ of course determine the same wall, but no other
numbers do.  Also, there are only finitely many walls, since for a
given $r^+$ and $\{ \ell_{j^+_i} \}$, the affine hyperplane defined by
the above equation only intersects $W$ for finitely many $d^+$.  The
connected components of the complement of the walls are the chambers;
for purely numerical reasons, the semistability condition is constant
in the interior of a chamber.

Now, as in the set-up of \S3, suppose $(\ell_j)$ lies on a single
wall in $W$, and choose $(\ell_j^+)$ and $(\ell_j^-)$ in the adjacent
chambers such that the line segment connecting them crosses a wall
only at $(\ell_j)$.  Then $x \in X$ belongs to $X^0$ if and only if it
represents a parabolic bundle which splits as $E_+ \oplus E_-$, where
$E_\pm$ are $(\ell_{j^\pm_i})$-stable parabolic bundles.  This is
because, to be in $X^{ss}(0) \sans X^{ss}(+)$, a parabolic bundle must
have a semistable parabolic subbundle $E_+$ of rank $r^+$, degree
$d^+$, and weights $\ell_{j^+_i}$.  Indeed, $E_+$ must be stable, for
since $(\ell_j)$ lies on only one wall, $E_+$ can have no parabolic
subbundle of the same slope.  For the same reason $E / E_+$ must be
stable.  On the other hand, to be in $X^{ss}(0) \sans X^{ss}(-)$, $E$
must have an stable parabolic subbundle $E_-$ of rank $r^-$, degree
$d^-$, and weights $\ell_{j^-_i}$.  Since all the weights are
distinct, $E_-$ cannot be isomorphic to $E_+$; so there is a nonzero
map $E_- \to E / E_+$.  By \cite[III Prop.\ 9(c)]{s}, this map must be
an isomorphism, so $E$ splits as $E_+ \oplus E_-$.

On the other hand, if $E_+$ and $E_-$ are any stable parabolic bundles
with rank, degree, and weights as above, then $E_+ \oplus E_-$ is
certainly represented in $X^0$. Hence $X^0 \mod G(0) = M(\ell_{j_i^+})
\times M(\ell_{j_i^-})$, the product of two smaller moduli spaces.

It is now easy to verify the hypotheses of our strongest result
\re{1x}.  First, $X$ is smooth on $X^{ss}(\ell_j)$, hence at $X^0$.
Second, for any $x \in X^0$, the stabilizer $G_x$ is the subgroup
isomorphic to $\kst$ acting on $H^0(E_+)$ with weight $\chi^- / c$ and
on $H^0(E_-)$ with weight $-\chi^+ / c$, where $c$ is the greatest
common divisor $(\chi^+, \chi_-)$.  This is because any $g \in
\GL{\chi}$ stabilizing a point in $X^0$ induces an automorphism of
$\be_x$, and vice versa, so for $x \in X^0$ there is an isomorphism
$\GL{\chi}_x \cong \Par \Aut (E_+ \oplus E_-) = \kst \times \kst$; but
only the automorphisms acting with the weights above correspond to
special linear $g$.  Third, if $L_j$ is the ample generator of $\Pic
\Gr(\chi - j, \chi)$, then for $x \in X^0$, $G_x \cong \kst$ acts on
$(L_j)_x$ with weight $(n_j^+ \chi^- - n_j^- \chi^+) / c$, where
$n_j^\pm$ are the number of $j_i^\pm$ less than or equal to $j$.  But
since $\chi^+ / c$ and $\chi^- / c$ are coprime, so are these weights
for some two values of $j$.  Bhosle gives a formula for the
linearization on $Z \times \Gr$ determined by $(\ell_j)$ in terms of
the $L_i$; an easy argument using this formula shows that
$(\ell_j^\pm)$ can be chosen within their chambers so that $G_x$ acts
with coprime weights on the corresponding linearizations.

Because all semistable parabolic bundles are represented by points in
$X$, and because semistability is an open condition, the universal
family of parabolic bundles is a versal family near any point $x \in
X^0$.  Moreover, two points in $X$ represent the same parabolic bundle
if and only if they are in the same orbit.  The normal bundle $N_{G
\cdot x / X}$ to an orbit is therefore exactly the deformation space
of the parabolic bundle.  For $\be_x = E_+ \oplus E_-$ as above this
is $H^1(X, \Par \End (E_+ \oplus E_-))$.  The stabilizer $G_x$ acts
with weight $r^\pm / (r^+, r^-)$ on $E_\pm$, so $N^0 = H^1(X, \Par
\End E_+) \oplus H^1(X, \Par \End E_-)$, and $N^\pm = H^1(X, \Par \Hom
(E_\mp, E_\pm))$.  Moreover, every element in $N^\pm$ is acted on with
weight exactly $\pm(\chi^+ + \chi_-)/ c$, so $N^\pm$ descend to vector
bundles $W^\pm$ over $M(\ell_{j_i^+}) \times M(\ell_{j_i^-})$.
Indeed, if $\be_\pm \to M(\ell_{j_i^\pm}) \times C$ are universal
bundles, then $W^\pm = (R^1 \pi) \Par \Hom(\be_\mp, \be_\pm)$.
Theorem \re{1x} then states that $\Pj W^\pm$ are the exceptional loci
of the morphisms $M(\ell_j^\pm) \to M(\ell_j)$.  This is the result of
Boden and Hu \cite{bhu}.  Moreover, \re{1x} asserts that the blow-ups
of $M(\ell_j^\pm)$ at $\Pj W^\pm$, and of $M(\ell_j)$ at
$M(\ell_{j_i^+}) \times M(\ell_{j_i^-})$, are all naturally isomorphic
to the fibred product $M(\ell_j^-) \times_{M(\ell_j)} M(\ell_j^+)$.

With the obvious modifications, the same techniques and results go
through for bundles with parabolic structure at several marked points,
or with degenerate flags.

\bit{Bradlow pairs}

The moduli spaces of Bradlow pairs on our curve $C$ can be studied in
the same way. The role of the weights will be played by a positive
parameter $\si \in \Q$.

A {\em Bradlow pair} is a pair $(E, \phi)$ consisting of a vector
bundle $E$ over $C$ and a nonzero section $\phi \in H^0(X,E)$.  We
refer to \cite{brad,bd,t1} for basic definitions and results on Bradlow
pairs.  As in \cite{t1}, we confine ourselves to the study of rank 2
pairs.  In this case a Bradlow pair of degree $d$ is $\si$-{\em
semistable} if for all line bundles $L \subset E$,
$$\begin{array}{cl}
\deg L \leq d/2 - \si & \mbox{if $\phi \in H^{0}(L)$ and} \\
\deg L \leq d/2 + \si & \mbox{if $\phi \not\in H^{0}(L)$.}
\end{array}$$
It is  $\si$-{\em stable} if both inequalities are strict.

The moduli spaces $B_d(\si)$ of $\si$-semistable rank 2 pairs were
constructed in \cite{t1}.  In that paper, the determinant was fixed,
but to parallel the discussion of parabolic bundles above we shall now
allow arbitrary determinant.  With that modification, the construction
goes as follows.

It suffices to construct $B_d(\si)$ for $d$ sufficiently large.  This
is because, for any effective divisor $D$, $B_d(\si)$ will be embedded
in $B_{d+ 2|D|}(\si)$ as the locus where $\phi$ vanishes on $D$.  So
assume that $d/2 - \si > 2g-2$, and let $\chi = d + 2(1-g)$.  Let
$\Quot$, $R$ and $Z$ be as in \S7 above, and let $G = \SL{\chi}$ act
diagonally on $R \times \Pj^{\chi-1}$.  The hypothesis $d/2 - \si >
2g-2$ implies that every $\si$-semistable pair is represented by a
point in $R \times \Pj^{\chi-1}$.  Let $(Z \times
\Pj^{\chi-1})^{ss}(\si)$ denote the semistable set with respect to the
fractional linearization $\co(\chi + 2 \si, 4 \si)$, and let $(R
\times \Pj^{\chi-1})^{ss}(\si)$ denote the $\si$-semistable set in the
sense of the definition above.  Then the natural $G$-morphism $T
\times 1: R \times \Pj^{\chi-1} \to Z \times \Pj^{\chi-1}$ satisfies
$T^{-1} (Z \times \Pj^{\chi-1})^{ss}(\si) = (R \times
\Pj^{\chi-1})^{ss}(\si)$.  Moreover, by \re{2u}, it is an embedding.
So if $X$ denotes the Zariski closure of its image, then for reasons
like those given in \S7, the moduli space $B_d(\si)$ is the geometric
invariant theory quotient $X \mod G(\si)$, where $\si$ denotes the
fractional linearization $T^* \co(\chi + 2 \si, 4 \si)$.

So again our main results apply.  The stability condition only changes
for $\si \in d/2 + \Z$, so these points are the walls.  Fix one such
$\si$.  Then $x \in X$ belongs to $X^0$ if and only if it represents a
pair which splits as $L \oplus M$, where $\deg L = d/2 - \si$ and
$\phi \in H^0(L)$.  Indeed, a subbundle $L$ of degree $d/2 - \si$ is
needed to violate the first semistability condition to the right of
the wall, and a subbundle $M$ of degree $d/2 + \si$ is needed to
violate the second semistability condition to the left.  But since
$\deg M = \deg E/L$, the map $M \to E/L$ is an isomorphism, so $E$ is
split.  On the other hand, for any pair $(L, \phi)$ with $\deg L = d/2
- \si$ and $\phi \in H^0(L) \sans 0$, and for any line bundle $M$ with
$\deg M = d/2 + \si$, certainly $(L \oplus M, \phi \oplus 0)$ is
represented in $X^0$.  Hence $X^0 \mod G(0) = S^iC \times \Pic^{d-i}
C$, where $i = d/2 - \si$ and $S^iC$ is the $i$th symmetric product.

It is now easy to verify the hypotheses of our strongest result
\re{1x}.  First, $X$ is smooth at $X^{ss}(0)$, hence at $X^0$.
Second, for any $x \in X^0$, the stabilizer $G_x$ is the subgroup
isomorphic to $\kst$ acting on $H^0(L)$ with weight $\chi(M)/c$ and on
$H^0(M)$ with weight $-\chi(L)/c$, where $c$ is the greatest common
divisor $(\chi(L), \chi(M))$.  This is because any $g \in \GL{\chi}$
stabilizing a point in $X^0$ induces an automorphism of the
corresponding pair, and vice versa, but only the automorphisms acting
with the weights above correspond to special linear $g$.  Third, this
stabilizer $\kst$ acts on $\co(1,0)_x$ with weight $(\chi(M) -
\chi(L))/c$, and on $\co(0,1)_x$ with weight $\chi(M)/c$.  These are
coprime, so linearizations with coprime weights can be chosen within
the chambers adjacent to $\si$.

As in \S7, the normal bundle to an orbit is exactly the deformation
space of the Bradlow pair.  For any pair $(E, \phi)$, this is the
hypercohomology group $\Hyp^1(\End E \stackrel{\phi}{\to} E)$.  (See
\cite{bd} or \cite[(2.1)]{t1}; the slightly different formula in
\cite{t1} arises because the determinant is fixed.)  More naturally,
the term $E$ in the complex is actually $E \otimes \co$, where $\co$
is the dual of the subsheaf of $E$ generated by $\phi$; it is acted on
accordingly by $G$.  For a pair $(L \oplus M, \phi \oplus 0)$
represented in $X^0$, this splits as
$$\Hyp^1(\co \oplus \co \stackrel{\phi \oplus 0}{\lrow} L)
\oplus \Hyp^1(LM^{-1} \to 0)
\oplus \Hyp^1(ML^{-1} \stackrel{\phi}{\to} M).$$

These are acted on by the stabilizer $\kst$ with weights $0$ and
$\pm(\chi(M) + \chi(L))/c$ respectively, so they are exactly $N^0$,
$N^+$, and $N^-$.  The expressions for $N^\pm$ can be simplified:
$\Hyp^1(LM^{-1} \to 0)$ is just $H^1(LM^{-1})$, and if $D$ is the
divisor of zeroes of $\phi$, the long exact sequence of $$ 0 \lrow
ML^{-1} \lrow M \lrow \co_D \otimes M \lrow 0 $$ implies that
$\Hyp^1(ML^{-1} \stackrel{\phi}{\to} M)$ is just $H^0(\co_D \otimes
M)$.  Since every element in $N^\pm$ is acted on with weight exactly
$\pm (\chi(M) + \chi(L))/c$, $N^\pm$ descend to vector bundles $W^\pm$
over $S^iC \times \Pic^{d-i}C$.  Indeed, if ${\bf M} \to \Pic^{d-i}C
\times C$ is a Poincar\'e bundle and $\Delta \subset S^iC \times C$ is
the universal divisor, then $W^+ = (R^1 \pi) {\bf M}^{-1}(\Delta)$ and
$W^- = (R^0 \pi) \co_\Delta \otimes {\bf M}$.

Theorem \re{1x} then states that $\Pj W^\pm$ are the exceptional loci
of the morphisms $B_d(\si \pm \half) \to B_d(\si)$, and that the
blow-ups of $B_d(\si \pm \half)$ at $\Pj W^\pm$, and of $B_d(\si)$ at
$S^iC \times \Pic^{d-i} C$, are all isomorphic to the fibred product
$B_d(\si - \half) \times_{B_d(\si)} B_d(\si + \half)$. This includes
the main result (3.18) of \cite{t1}; to recover the $W^\pm$ obtained
there for a fixed determinant line bundle $\La$, just substitute ${\bf
M} = \La(-\Delta)$.

Notice that since the construction only works for $d$ large, the
result has so far only been proved in that case.  For general $d$,
choose disjoint divisors $D$, $D'$ of degree $|D| = |D'|$ such that $d
+ 2|D|$ is large enough.  Then $D$ and $D'$ determine two different
embeddings $B_{d+2|D|}(\si) \to B_{d+4|D|}(\si)$ whose images
intersect in $B_d(\si)$.  The result for $B_d(\si)$ then follows
readily from the result for $B_{d+2|D|}(\si)$ and $B_{d+4|D|}(\si)$. \bl

A similar argument proves the analogous result of Bertram,
Daskalopoulos and Wentworth \cite{bdw} on Bradlow $n$-pairs.  These
are pairs $(E, \phi)$, where $E$ is as before, but $\phi$ is now a
nonzero element of $H^0(E \otimes \co^n)$.  The stability condition is
just like that for ordinary Bradlow pairs, except that the two cases
are $\phi \in H^0(L \otimes \co^n)$ and $\phi \notin H^0(L \otimes
\co^n)$.  There is no geometric invariant theory construction in the
literature of moduli spaces of $n$-pairs, but the construction of
\cite{t1} for 1-pairs generalizes in the obvious way; for example,
$\Pj^{\chi-1}$ gets replaced by $\Pj^{n\chi-1}$.

Again the stability condition only changes for $\si \in d/2 + \Z$, so
fix one such $\si$.  Assume that $d/2 - \si > 2g-2$, so that the
moduli space can be constructed directly as a geometric invariant
theory quotient.  Then $x \in X^0$ if and only if it splits as $L
\oplus M$, where $\deg L = i = d/2 - \si$ and $\phi \in H^0(L \otimes
\co^n)$.  Hence $X^0 \mod G(0) = \times^n_{\Pic^i C} S^i C
\times \Pic^{d-i}C$, where $\times^n_{\Pic^i C}$ denotes the $n$-fold
fibred product over $\Pic^i C$.  The hypotheses of \re{1x} are
verified exactly as before.  The normal bundle to an orbit is again
the deformation space, but this is now $\Hyp^1(\End E
\stackrel{\phi}{\to} E \otimes \co^n)$; for a pair $(L \oplus M, \phi
\oplus 0)$ represented in $X^0$, this splits as
$$\Hyp^1(\co\oplus \co \stackrel{\phi\oplus 0}{\lrow} L \otimes \co^n)
\oplus \Hyp^1(LM^{-1} \to 0)
\oplus \Hyp^1(ML^{-1} \stackrel{\phi}{\to} M \otimes \co^n).$$
Again these are exactly $N^0$, $N^+$, and $N^-$.  The expression for
$N^+$ is simply $H^1(LM^{-1})$, but the expression for $N^-$
cannot be simplified very much. If $F$ is defined to be the cokernel
of the sheaf injection $ML^{-1} \to M \otimes \co^n$ induced by
$\phi$, then $N^- = H^0(F)$, but this is not very helpful as $F$ may
not be locally free.  In any case, $N^\pm$ descend as before.  If
$\Delta_j \subset \times^n_{\Pic^i C} S^i C \times C$ is the pullback
from the $j$th factor of the universal divisor $\Delta \subset S^iC
\times C$, then there is a universal map $\co \to \oplus_j
\co(\Delta_j)$ of bundles on $\times^n_{\Pic^i C} S^i C \times C$, and
hence a sheaf injection ${\bf ML}^{-1} \to {\bf ML}^{-1} \otimes
\oplus_j \co(\Delta_j)$ of bundles on $\times^n_{\Pic^i C} S^i C
\times \Pic^{d-i} C \times C$, where ${\bf M}$ and ${\bf L}$ are
Poincar\'e bundles on $\Pic^i C \times C$ and $\Pic^{d-i}C \times C$
respectively.  Let ${\bf F}$ be the cokernel of this injection; then
$W^- = (R^0\pi) {\bf F}$.  As before, $W^+ = (R^1 \pi){\bf LM}^{-1}$.

Theorem \re{1x} then gives a result precisely analogous to the one
stated above for 1-pairs.  However, the argument passing from large
$d$ to general $d$ no longer works, so this result is only valid for
$i = d/2 - \si > 2g-2$.  Indeed, for $i$ smaller than this, the fibred
product $\times^n_{\Pic^i C} S^i C$, and hence $X^0 \mod G(0)$, have
bad singularities.  Nevertheless, Bertram et al.\ succeed in using
this result to compute certain Gromov invariants of Grassmannians.

{\em Acknowledgements.}  I am grateful to Aaron Bertram, Igor
Dolgachev, David Eisenbud, Jack Evans, Antonella Grassi, Stuart
Jarvis, J\'anos Koll\'ar, Kenji Matsuki, Rahul Pandharipande, David
Reed, Eve Simms, and especially Frances Kirwan, Peter Kronheimer, and
Miles Reid for very helpful conversations and advice.

\end{document}